\documentclass{amsart}
\usepackage[
  marginparwidth=5mm,     
  marginparsep=1mm       
  ]{geometry}
\usepackage{graphicx} 
\usepackage{fullpage}

\usepackage[utf8]{inputenc}
\usepackage[T1]{fontenc}
\usepackage{lmodern}
\usepackage{amssymb}
\usepackage{multicol}
\usepackage{longtable}

\usepackage{newfloat}
\DeclareFloatingEnvironment{algorithm}

\usepackage{amsthm}
\newtheorem{theorem}{Theorem}
\newtheorem{lemma}[theorem]{Lemma}
\newtheorem{prop}[theorem]{Proposition}
\newtheorem{corollary}[theorem]{Corollary}
\newtheorem{defn}[theorem]{Definition}

\newcommand{\bN}{\mathbb N}
\newcommand{\bQ}{\mathbb Q}
\newcommand{\bZ}{\mathbb Z}

\newcommand{\cH}{\mathcal H}

\newcommand{\serF}{\bar{F}}
\newcommand{\serf}{\bar{f}}

\newcommand{\scalp}[2]{\langle#1,#2\rangle}
\newcommand{\Scalp}[2]{\biggl\langle#1,#2\biggr\rangle}
\newcommand{\adj}[1]{#1^\dagger}
\newcommand{\twist}[1]{#1^\sharp}
\newcommand{\ann}{\operatorname{ann}}
\newcommand{\red}{\operatorname{red}}
\newcommand{\lt}{\operatorname{lt}}

\newcommand{\vardbtilde}[1]{\tilde{\raisebox{0pt}[0.85\height]{$\tilde{#1}$}}}

\usepackage[shortlabels]{enumitem}

\usepackage{xcolor,soul}

\usepackage{url}
\usepackage[plainpages=false,pdfpagelabels,colorlinks=true,allcolors=blue!80!black,hypertexnames=false,pdftex=true,unicode]{hyperref}

\usepackage{amsaddr}
\title{Differential equations satisfied by generating functions\\
  of 5-, 6-, and 7-regular labelled graphs:\\
  a reduction-based approach}
\author{Frédéric Chyzak}
\address{Inria, France}
\email{frederic.chyzak@inria.fr}
\author{Marni Mishna}
\address{Simon Fraser University, Canada}
\email{mmishna@sfu.ca}
\date{\today}

\begin{document}
\begin{abstract}
  By a classic result of Gessel, the exponential generating functions for
  $k$-regular graphs are D-finite. Using Gr\"obner bases in Weyl
  algebras, we compute
  the linear differential equations satisfied by the generating
  function for 5-, 6-, and 7- regular graphs. The method is sufficiently
  robust to consider variants such as graphs with multiple edges,
  loops, and graphs whose degrees are limited to fixed sets of values.

\bigskip
\noindent\emph{Keywords:} regular graph, enumeration, Weyl algebra, reduction-based integration.

\noindent\emph{2020 MSC:} 05C30, 12H05.
\end{abstract}
\maketitle

\section{Introduction}
\label{sec:introduction}
\subsection{A short history of $k$-regular graph enumeration}
A graph is said to be \emph{regular} if every vertex is incident to
the same number of edges, that is, each vertex has the same degree. If
that degree is $k$, we call the graph $k$-regular. One of the earliest
graph enumeration problems considered was the number of non-isomorphic unlabelled $k$-regular graphs on $n$ vertices.
It is a relatively attainable problem for many reasons, including the fact that the number of edges is fixed in these graphs,
which yields a significant simplification.  For example, according to Gropp~\cite{gropp_enumeration_1992}, Jan de Vries determined the number of non-isomorphic cubic (3-regular) graphs up to 10 vertices, and shared them in a letter to Vittorio Martinetti, which was eventually published in a journal in 1891. The proofs were descriptions of the graphs. Here we consider the slightly easier problem of labelled graphs, specifically the number of labelled $k$-regular graphs on~$n$ vertices, which we denote by~$r^{(k)}_n$.

In the labelled case, the work of Read in the 1950s established enumeration formulas using the cycle index series, a relatively new machinery at the time. He gives a compact, structural equation in~\cite[Eq.~5.11]{read_enumeration_1959} that is not immediately suitable for enumeration purposes for $k>3$. He notes,
 \begin{quote}\emph{
“It may readily be seen that to evaluate the above expressions in particular cases may involve an inordinate amount of computation.”}
\end{quote}
For $k=3$, the equation is sufficiently manageable to give rise to a nice asymptotic formula.

One can distill from his work a formula in terms of coefficient
extraction of a multivariable polynomial. This is the starting point of most modern approaches as it is easy to interpret, and there are numerous possibilities for analysis.
Using the notation of square brackets to isolate the coefficient of
the indicated term in a series expansion of the product we can
write
\begin{equation}\label{eq:extraction}
r^{(k)}_n = [x_1^kx_2^k\dots x_n^k] \prod_{1\leq i<j\leq n} (1+x_ix_j).
\end{equation}
The multiplication accounts for all possibilities of an edge $\{i,j\}$ to be in the graph or not. The coefficient of the indicated monomial is the number of graphs that have vertices $1$ to $n$, such that each vertex is incident to exactly $k$ other vertices: this is precisely $r^{(k)}_n$.

There are a variety of strategies to consider for this construction,
and other direct arguments in the service of asymptotic enumeration.
The problem is well studied, remarkably even for problems with $k$ a function of~$n$.
Wormald's 2018 ICM survey has many details on the
state of asymptotic enumeration of regular graphs and related
objects~\cite{wormald_asymptotic_2018}.

In that 2018 survey, Wormald notes that no new exact enumeration
results have appeared since the recurrences for 4-regular graphs
published in the early 1980s. The entry point of the present article
is also Eq.~\eqref{eq:extraction}, but we follow a different lineage
to contribute fixed-length linear recurrence formulas to count 5-,~6-, and~7-regular graphs, ending the drought.

The fact that there are recurrences to find at all is related to a
question of Stanley~\cite{stanley_differentiably_1980} in his
foundational article on P-recursive sequences. The existence of a
recurrence is equivalent to asking whether or not the exponential
generating function for $r^{(k)}_n$, defined as
$R^{(k)}(t):= \sum_{n\geq 0}r^{(k)}_n\frac{t^n}{n!}$, is
\emph{D-finite}.  In other words, does $R^{(k)}(t)$ satisfy a linear
differential equation with polynomial coefficients? Read had already
given a recurrence for 3-regular graphs in his PhD thesis~\cite{read_some_1959}, and Read
and Wormald used a combinatorial analysis to produce recurrences for
4-regular graphs~\cite{read_number_1980}.  Goulden, Jackson and
Reilly~\cite{goulden_hammond_1983} were also able to determine
explicit linear differential equations satisfied by $R^{(3)}(t)$ and
$R^{(4)}$ using tools that dated back to MacMahon at the turn of the
20th century, called Hammond operators. But, they noted
that\footnote{In our notation, $p = k$ and $y_p=R^{(k)}$.}
\begin{quote}\emph{
    ``\dots the H-series theorem enables us to write down the system of partial differential equations for the H-series for arbitrary $p$ without difficulty. However, the reduction of this system to a single ordinary differential equation in $y_p$ is a technical task which we are unable to carry out for the general case.''}
\end{quote}
Their work fuelled speculation that $R^{(k)}$ should be D-finite for all~$k$.
Gessel compared their approach to his own method by the scalar product of symmetric functions and algebraic substitutions~\cite{gessel_enumerative_1987}:
\begin{quote}\emph{
``\dots Hammond operators are undesirable for two reasons. First, they disguise the symmetry of the scalar product. Second, they can be represented as differential operators. Although this might seem like an advantage, it seems to be of little use, but misleads by directing attention in the wrong direction.''}
\end{quote}
Instead of working with differential equations, he recast the
extraction in terms of symmetric functions, and used algebraic
arguments to establish that indeed $R^{(k)}(t)$ is D-finite for all
$k$. His framework is sufficiently simple and robust that it can be
used to establish the D-finiteness of many related regular graph and
hypergraph cases.  Gessel was able to advance on the general case
thanks to concurrent work on multivariable P-recursiveness of
Lipshitz~\cite{Lipshitz-1989-DFP}. The work of Lipshitz was not
sufficiently straightforward to convert into an algorithm or even make
computation effective beyond $k=2$. It was over a decade before the
computer algebra implementations using differential operators caught
up to his theoretical results.
In 2005 Chyzak, Mishna and Salvy~\cite{chyzak_effective_2005} made both the Hammond method
and the Gessel strategy effective for any~$k$ using Gröbner bases for
D-modules and non-commutative polynomial elimination, in a sort of
variant of Creative Telescoping, a method for symbolic integration.
The implementation quickly found differential
equations up to, and including, 4-regular objects.  The
growth of data in the skew polynomial elimination involved in the
5-regular graph case requires computational resources that even today
are insufficient to have the algorithm terminate.  However, in the
intervening 20 years, there have been remarkable improvements and insights to Creative
Telescoping.
This lead us to an evolved algorithm that terminates also
in practice, and indeed we could find the linear differential
equations satisfied by $R^{(5)}(t)$, $R^{(6)}(t)$, and $R^{(7)}(t)$.
Our present approach can be applied to find the differential equations
satisfied by the other graph, hypergraph and graph-like classes
for higher degrees of regularity than were previously 
obtained~\cite{mishna_automatic_2007,mishna_regularity_2018}.

The following theorem is the main result of this article.
It appears below, rephrased, as Corollary~\ref{cor:table-is-correct}(\ref{item:lde-is-correct}).

\begin{theorem}\label{thm:ldes}
For each graph model in Table~\ref{tab:results}, there exists a
known linear differential equation with polynomial coefficients
satisfied by the exponential generating function,
with explicit order given by column~$\partial_t$ of the table,
and maximum coefficient degree given by column~$t$ of the table.

The graph models in the table include those with
\begin{enumerate}
\item only simple edges permitted (denoted `se' in the table) or  multiple edges allowed (denoted `me');
\item loops forbidden (`ll'), loops allowed and contributing~2 to vertex degrees (`la'), or loops allowed and contributing~1 to vertex degrees (`lh');
\item degrees restricted to some finite set, including: $\{k\}$
  for~$2\leq k\leq 7$,  $\{1,\dots,k\}$ for~$2\leq k\leq 6$, and
  $\{k,\ell\}$ for~$2\leq k < \ell\leq 6$, among others.
\end{enumerate}

In particular, the order of the linear differential equation for
simple, loopless $k$-regular graphs (coded `se' and `ll') is
summarized in the table:
\begin{center}
\begin{tabular}{l|llllll}
$k$&2&3&4&5&6&7\\\hline
order&1&2&2&6&6&20
\end{tabular}.
\end{center}
\end{theorem}

\subsection{The scalar product%
\protect\footnote{We follow the usual terminology of a “scalar product”
in combinatorics,
although the presence of a formal indeterminate~$t$
would require to speak more properly of a “pairing”.} %
of symmetric functions}
\label{sec:scalarproduct}

The coefficient extraction in Eq.~\eqref{eq:extraction} can be placed
into an infinite product, symmetric in all variables, which can be
readily encoded in terms of symmetric functions. The set up of
Gessel~\cite{gessel_symmetric_1990} uses the scalar product in the
ring of symmetric functions to model the extraction.  Describing the
method requires a small detour through symmetric function terminology
and basics. There are many excellent introductions. We highlight some
notation, but refer readers to
Stanley~\cite[Chapter 7]{stanley_enumerative_1999} for
details.

We say $\lambda=(\lambda_1,\lambda_2\dots, \lambda_q)$ such that $\sum_{i=1}^q\lambda_i=n$ and $\lambda_i\geq \lambda_{i+1}$ is a partition of~$n$ into $q$~parts,  and write
$\lambda\vdash n$ to indicate that $\lambda$ is a partition of~$n$. The monomial symmetric
function is defined as $m_\lambda:=\sum_{\alpha\sim\lambda}
x^\alpha$
where $\alpha\sim\lambda$ if the non-zero entries of $\alpha$ are a rearrangement
of the parts of $\lambda$. Using $m_\lambda$ we can
describe the complete homogeneous symmetric function
$h_n:=\sum_{\lambda \vdash n} m_\lambda$
and the power-sum symmetric function
$p_n:= m_{(n)}=x_1^n+x_2^n+\dots$.
Products are denoted respectively
$h_{n_1\,n_2\dots\,n_{\ell}}:=h_{n_1}h_{n_2}\dots h_{n_\ell}$ and
$p_{n_1\,n_2\dots\,n_{\ell}}:=p_{n_1}p_{n_2}\dots p_{n_\ell}$. The
vector space of symmetric functions of order~$n$ has numerous bases,
including $\{m_\lambda\mid \lambda\vdash n\}$,
$\{h_\lambda\mid \lambda\vdash n\}$ and  $\{p_\lambda \mid \lambda\vdash n\}$.
For any $\lambda\vdash n$,  $z_\lambda$ denotes the number
\begin{equation}\label{eq:z-lambda}
z_\lambda := 1^{r_1} r_1!\, 2^{r_2} r_2!\, \dots n^{r_n} r_n!
\end{equation}
provided $\lambda$ has $r_1$ ones, $r_2$ twos, etc.
We set $\delta_{\lambda=\nu}$ to~1 if $\lambda=\nu$ is true
and to~0 otherwise.
The \emph{scalar product of symmetric functions} is classically defined by
\begin{equation}
\scalp{p_\lambda}{p_\nu} := \delta_{\lambda=\nu} z_\lambda ,
\quad \text{from which we deduce}\quad
\scalp{m_\lambda}{h_\nu} =\delta_{\lambda=\nu} .
\end{equation}
The connection to the graph enumeration problem is as follows.
We can
extract the coefficient of a particular monomial in a symmetric
function with a judiciously chosen scalar product. Write
$\serF:=\prod_{i<j} (1+x_ix_j) $ and consider an example. This
product is fundamental in the study of symmetric functions,
particularly its expression in the various bases. Now, since
$r^{(3)}_4=[x_1^3x_2^3x_3^3x_4^3] \serF$,  to actually compute
this write $\serF$ as a
sum of monomial symmetric functions, and determine the coefficient of
$m_{3,3,3,3}$ (which is the only basis element to contain the term
$x_1^3x_2^3x_3^3x_4^3$). As the monomial and complete homogenous bases
are orthogonal under the usual scalar product of symmetric functions, this coefficient is precisely
$\scalp{\serF}{h_{3,3,3,3}} = \scalp{\serF}{h_3^4}$.

From the formulas $\log(1+u) = \sum_{k\geq1} (-1)^{k+1}u^k/k$
and $2\sum_{i<j}x_i^kx_j^k = \sum_{i,j}x_i^kx_j^k - \sum_ix_i^{2k}$
it follows
\begin{equation}\label{eqn:Gequiv}
\serF 
= \exp\Biggl(\sum_{i<j} \log(1+x_ix_j)\Biggr)
= \exp\Biggl(\sum_{i<j} \sum_{k\geq1} (-1)^{k+1}\frac{x_i^kx_j^k}{k}\Biggr)
= \exp\Biggl(\sum_{k\geq1} (-1)^{k+1}\frac{p_k^2-p_{2k}}{2k}\Biggr) .
\end{equation}

 Henceforth we will only work with the power-sum basis,
 specifically, we work in a ring generated by~$t$ and a finite number
 of the the $p_i$ variables. To continue the example, to determine
 $R^{(3)}(t)$ we first write
 $h_3=\frac{p_3}{3}+\frac{p_2p_1}{2}+\frac{p_1^3}{6}$, and thus obtain
 the following expression for the generating function:
 \begin{equation}\label{eq:R3-as-scalproduct-of-G}
   R^{(3)}(t)
   = \Scalp{\serF}{\sum_{n\geq0} h_3^n \frac{t^n}{n!}}
   = \Scalp{\exp\Biggl(\sum_{k\geq1}(-1)^{k+1}\frac{p_k^2-p_{2k}}{2k}\Biggr)}
           {\exp\left(\left(\frac{p_3}{3}+\frac{p_2p_1}{2}+\frac{p_1^3}{6}\right)t\right)} .
  \end{equation}
Since the second argument has only $p_1, p_2, p_3$, all terms with other
$p_i$ contribute 0:
\begin{equation}\label{eq:R3-as-scalproduct}
R^{(3)}(t) = \Scalp{\exp\left(\frac{p_1^2}{2}-\frac{p_2}{2}-\frac{p_2^2}{4}+\frac{p_3^2}{6}\right)}
                   {\exp\left(\left(\frac{p_3}{3}+\frac{p_2p_1}{2}+\frac{p_1^3}{6} \right)t\right)} .
\end{equation}

For future reference, we note the following formula,
which leads to generalizations
of Eqs.~\eqref{eq:R3-as-scalproduct-of-G} and~\eqref{eq:R3-as-scalproduct}:
 \begin{equation}\label{eq:Rk-as-scalproduct-of-G}
   R^{(k)}(t) =
     \Scalp{\serF}{\sum_{n\geq0} h_k^n \frac{t^n}{n!}} =
     \scalp{\serF}{\exp(h_k t)} .
 \end{equation}
This formula will be proven and extended in Lemma~\ref{lem:scalp-of-f-and-g}.

\subsection{Earlier computational approaches}
As we mentioned above, Gessel~\cite{gessel_symmetric_1990} proved
the existence of linear differential equations
for scalar products like Eq.~\eqref{eq:Rk-as-scalproduct-of-G},
and earlier work~\cite{chyzak_effective_2005}
proposed algorithms to compute them.
In there, for a given series~$S$ in the variables $p_1,\dots,p_k$
we consider the set, denoted~$\ann(S)$,
of all linear differential operators that annihilate~$S$.
The elements of~$\ann(S)$ are non-commutative polynomials
in the variables $p_1,\dots,p_k$
and in the corresponding derivatives $\partial_1,\dots,\partial_k$;
they possess a well-defined total degree in the $2k$ variables.
The set~$\ann(S)$ is closed under multiplication by any operator on the left
and is thus a left ideal.
As is customary in effective literature,
such a left ideal is best represented by a non-commutative analogue
of a Gröbner basis,
that is, by a finite set of non-commutative polynomials
that can algorithmically divide a given ideal element,
resulting into a uniquely defined remainder
that is zero if and only if the given polynomial is in the ideal.

Given a number~$k$, we henceforth write
$p = (p_1,\dots,p_k)$ and $\partial = (\partial_1,\dots,\partial_k)$.
Given a series~$F$ in~$p$ and a series~$G$ in~$(t,p)$,
the differential equations with respect to~$t$
satisfied by the scalar product~$\scalp{F}{G}$
are to be found as those elements free of~$(p,\partial)$
in the (vector space) sum
of the left ideal~$\ann(G)$ and of the right ideal~$\adj{\ann(F)}$
obtained by taking the adjoints of all elements in~$\ann(F)$~\cite{chyzak_effective_2005}
(see the definitions in Section~\ref{sec:approach}).
A first algorithm in~\cite{chyzak_effective_2005}, based on linear algebra,
consists:
\emph{(i)\/}~in fixing an integer~$d$;
\emph{(ii)\/}~in determining representatives of $\adj{\ann(F)}$ and~$\ann(G)$
for each possible leading monomial of total degree at most~$d$
with respect to~$(p,\partial,\partial_t)$;
\emph{(iii)\/}~and in using a non-commutative variant of Gaussian elimination
over~$\bQ(t)$
to eliminate~$(p,\partial)$,
repeating the whole process with a larger~$d$
if elimination results in no non-trivial output.
Because there are $\binom{d}{2k+1} = O(d^{2k+1})$ monomials of degree at most~$d$,
and almost as many representatives to determine for each ideal,
this process is very inefficient in practice.
A second algorithm in~\cite{chyzak_effective_2005} is tailored
to a certain form for the argument~$G$ in the scalar product:
if $G = \exp(h_kt)$, the theory of Hammond series,
as developed in~\cite{goulden_hammond_1983},
provides the formula
\[ \scalp{F}{\exp(h_kt_k)} = \cH(F)(0,\dots,0,t_k) , \] where
$\cH(F)(t_1,\dots,t_k)$ is a transform of~$F$ known as its Hammond
series.  A simple replacement of the~$p_i$ and the~$\partial_i$
in~$\ann(F)$ with suitable polynomials in $t_1,\dots,t_k$ and
corresponding derivatives~$\partial_{t_i}$ provides $\ann\bigl(\cH(F)\bigr)$.
The specialization of $t_1,\dots,t_{k-1}$ to~$0$ is then obtained by
restriction, an operation dual to integration.  One way to implement
it would have been to first eliminate the $k-1$ variables
$\partial_{t_1},\dots,\partial_{t_{k-1}}$, e.g., by a Gröbner basis
calculation, before setting all of the $k-1$ variables
$t_1,\dots,t_{k-1}$ to zero and taking a generator of the resulting
principal ideal in $\bQ(t_k)\langle\partial_{t_k}\rangle$.  But a
simultaneous elimination in this way leads to high degrees and is also
inefficient in practice.  More generally, in the 2000s, no good
algorithm was known for integration with respect to several variables
considered simultaneously, so one had to resort to iterated
integrations, one variable after the other.  Correspondingly, for
multiple restriction one had to perform specializations one variable
after the other, and this is what is proposed
in~\cite{chyzak_effective_2005}, in a way that is reminiscent of
elimination by successive resultants.  This approach, too, fails for
$k=5$: all steps are fast until the last elimination, which should
eliminate $\partial_{t_1}$ from two degree-9 polynomials in the four
variables $t_1,t_5,\partial_{t_1},\partial_{t_5}$, and this fails in
practice.

In both old approaches, the culprit is elimination in too many variables:
eliminating $2k$ variables between polynomials in $2k+1$ variables over~$\bQ(t)$
in the first approach;
eliminating $k-1$ variables between polynomials in $2k-1$ variables over~$\bQ(t_k)$
in the second approach.
The second is an improvement in that it reduces the number of variables,
and this is assisted by specializations to zero along the process.

A turning point in the theory of Creative Telescoping was the introduction
of reduction-based algorithms,
starting with the integration
of bivariate rational functions~\cite{BostanChenChyzakLi-2010-CCT} in 2010,
and followed by many articles in the literature.
Inspiration for the present work came
from a more recent reduction-based algorithm~\cite{BostanChyzakLairezSalvy-2018-GHR}
for the integration with regard to one variable~$p$
of general D-finite functions~$f(t,p)$,
leading to integrals parametrized by~$t$.
In a nutshell, reduction-based algorithms:
\emph{(i)\/}~set up a reduction process
that corresponds to simplifying a function to be integrated
modulo derivatives with respect to~$p$ of other functions,
in such a way that the resulting remainder lies
in a finite-dimensional vector space;
\emph{(ii)\/}~find a linear relation between the remainders
of successive higher-order derivatives with respect to the parameter~$t$
of the function to be integrated.
In situations where integrals of derivatives are zero,
the output linear relation reflects a differential equation in~$t$
of the parametrized integral.
Although the symmetric scalar product cannot be represented
as an integral of a D-finite function,
the method of~\cite{BostanChyzakLairezSalvy-2018-GHR}
can be adapted to the present situation, in a way
that the reduction with respect to the $k$ variables $p_1,\dots,p_k$
is possible simultaneously
and that most of the calculations involve polynomials in $k+1$ variables over~$\bQ(t)$.

\subsection{Contributions}

Besides presenting in Algorithm~\ref{algo:method} a method that adapts reduction-based
algorithms to a simultaneous reduction with respect to several integration variables,
our main contribution in the present work is to obtain differential equations
satisfied by various models of graphs with vertex degrees restricted to be
in a fixed subset of $\{1,\dots,7\}$
(see Theorem~\ref{thm:ldes} and Corollary~\ref{cor:table-is-correct}).
We cannot guarantee the termination of our method,
but any differential equation it outputs is correct,
as proven by Theorem~\ref{thm:proc-is-correct}.
In Table~\ref{tab:results}, we list
for a few dozens of models
the order of a differential equation satisfied by the counting generating function
and the order of a recurrence equation satisfied by its sequence of coefficients,
together with corresponding degrees of their coefficients.
All those equations are proven correct by the computer calculations (see Corollary~\ref{cor:table-is-correct}).
To the best of our knowledge, this is the first time differential equations
are presented for $R^{(5)}(t)$, $R^{(6)}(t)$ and~$R^{(7)}(t)$, or more generally
graphs where degrees 5, 6, or~7 are considered.

The recurrences we find are linear, with polynomial coefficients and
hence can be unravelled quickly to get data for graphs of high order.
For example, it take about 15~minutes
to determine the number of 7-regular graphs on 2000 vertices
from the ODE of order~20 that we found:
\[r^{(7)}_{2000} = 80680697\dots04296875 \approx 8.068069734\times10^{18572}.\]
It is even faster when the machine allows parallel processes.
More generally,
we are able to significantly increase the number of known terms
compared to the state of the art in the On-line Encyclopedia of
Integer Sequences (OEIS)~\cite{oeis} for sequences A338978 and~A339847, and we have contributed a new
sequence, A374842 counting 7-regular graphs.

The generated enumerative data, recurrences, differential equations
and Maple code implementing our strategy are all available at
\url{https://files.inria.fr/chyzak/kregs/}.

\section{Worked example: 4-regular graphs}

Before introducing our procedure in a systematic way
in Section~\ref{sec:approach},
we illustrate it with the class of 4-regular graphs,
allowing single edges and no loops.
(The case $k=3$ is too simple to demonstrate important points of our method.)
Specializing Eq.~\eqref{eq:Rk-as-scalproduct-of-G} to~$k=4$,
we consider the scalar product $\scalp{F}{G}$,
which represents~$R^{(4)}(t)$
when the exponential functions $F=\exp(f)$ and~$G=\exp(tg)$ are given by
\begin{equation*}
f := \frac{p_1^2}2 - \frac{p_2^2}4 + \frac{p_3^2}{6} - \frac{p_4^2}8 - \frac{p_2}2 + \frac{p_4}4 , \qquad
g := \frac{p_1^4}{24} + \frac{p_1^2 p_2}4 + \frac{p_2^2}8 + \frac{p_1 p_3}3 + \frac{p_4}4 .
\end{equation*}

\subsection{A reduction procedure}

We begin by explaining a procedure to normalize
expressions of the form $\scalp{F}{sG}$
for a polynomial $s \in \bQ (t)[p]$:
without changing the value of the scalar product,
the polynomial~$s$ will be replaced with an element in
$\bQ(t) + \bQ(t)p_1 + \bQ(t)p_2$.

From the definition of~$F$, we get that annihilating operators for~$F$ are
\begin{equation}\label{eq:Pi-for-4regs}
P_1 := \partial_1 - p_1 , \qquad
P_2 := 2\partial_2 + p_2 + 1 , \qquad
P_3 := 3\partial_3 - p_3 , \qquad
P_4 := 4\partial_4 + p_4 - 1 .
\end{equation}
In Section~\ref{sec:approach},
we will define two transformations on differential operators,
namely adjoints ($\adj{}$) and twists ($\twist{}$).
Applying them to Eq.~\eqref{eq:Pi-for-4regs},
we obtain
\begin{equation*}
\adj{P_1} := p_1 - \partial_1 , \qquad
\adj{P_2} := p_2 + 2\partial_2 + 1 , \qquad
\adj{P_3} := p_3 - 3\partial_3 , \qquad
\adj{P_4} := p_4 + 4\partial_4 - 1 ,
\end{equation*}
and
\begin{gather*}
\twist{P_1} := p_1 - \partial_1 - \frac{t}6(p_1^3+3p_1p_2+2p_3) , \qquad
\twist{P_2} := p_2 + 2\partial_2 + \frac{t}2(p_1^2+p_2) + 1 , \\
\twist{P_3} := p_3 - 3\partial_3 - t p_1 , \qquad
\twist{P_4} := p_4 + 4\partial_4 + t - 1.
\end{gather*}
We will prove in Section~\ref{sec:approach}
that $\scalp{F}{(\twist{P_j}\cdot \bar s) \, G}$ is zero
for any $\bar s \in \bQ(t)[p]$ and any~$j$,
motivating that we will try to adjust~$s$
by a linear combination of polynomials of the form $\twist{P_j}\cdot \bar s$.

In order to determine how to do so more precisely,
observe first that for any monomial~$p^\alpha$,
\begin{gather*}
\twist{P_1}\cdot p^\alpha = -\frac{t}6 p_1^{\alpha_1+3}p_2^{\alpha_2}p_3^{\alpha_3}p_4^{\alpha_4} + \dotsb , \qquad
\twist{P_2}\cdot p^\alpha = \frac{t}2 p_1^{\alpha_1+2}p_2^{\alpha_2}p_3^{\alpha_3}p_4^{\alpha_4} + \dotsb , \\
\twist{P_3}\cdot p^\alpha = -t p_1^{\alpha_1+1}p_2^{\alpha_2}p_3^{\alpha_3}p_4^{\alpha_4} + p_1^{\alpha_1}p_2^{\alpha_2}p_3^{\alpha_3+1}p_4^{\alpha_4} + \dotsb , \qquad
\twist{P_4}\cdot p^\alpha = p_1^{\alpha_1}p_2^{\alpha_2}p_3^{\alpha_3}p_4^{\alpha_4+1} + \dotsb ,
\end{gather*}
where in each case, the dots represent a polynomial with lower total degree.
We will base our calculation on these forms.
Consider for example any monomial ordering for which $p_4$~is
lexicographically higher than all other variables.
Given a polynomial~$s \in \bQ(t)[p]$
with leading term~$c p^\beta$ for~$\beta_4\geq1$,
the choice $\alpha = \beta-(0,0,0,1)$ ensures
that $s - \twist{P_4}\cdot(c p^\alpha)$ has a leading monomial
less than~$p^\beta$.
As a consequence,
$s$~can be reduced by a series of like transformations
to a polynomial $s - \twist{P_4}\cdot \bar s$
that does not involve~$p_4$:
here $\bar s$~is a polynomial that adds up all the $c p^\alpha$
observed during the reduction process.
In other words, one can eliminate~$p_4$ from~$s$.
One can similarly use~$\twist{P_3}$ to reduce
the degree with respect to~$p_3$:
this essentially introduces~$p_1$ as a replacement of~$p_3$,
but one can eliminate~$p_3$ as well.
By continuing with transformations based on~$\twist{P_2}$,
which do not reintroduce either $p_3$ or~$p_4$,
one could hope to eliminate~$p_1$ as well (after $p_3$ and~$p_4$) from~$s$.
It turns out
that one cannot fully eliminate~$p_1$,
but that degrees with respect to~$p_1$ can be reduced
down to at most~$1$.
On the other hand, it is not immediately evident
that degrees with respect to~$p_2$ can be kept under control.

To explain how controling~$p_2$ can be done,
we continue our informal presentation
by recombining the~$\twist{P_i}$
in the following way
into elements of the right ideal they generate:
\begin{align*}
\twist{P_1} + \twist{P_3}\frac{t}3 &= -\frac{t}6p_1^3 - \frac{t}2p_1p_2 + \left(1-\frac{t^2}3\right)p_1 - \partial_1 - t\partial_3 , \\
\twist{P_2} &= \frac{t}2p_1^2 + \left(1+\frac{t}2\right)p_2 + 1 + 2\partial_2 , \\
\tilde P_5 :=
\twist{P_1} + \twist{P_3}\frac{t}3 + \twist{P_2}\frac{p_1}3 &= \frac{1-t}3p_1p_2 + \frac{4-t^2}3p_1 + \frac23p_1\partial_2 - \partial_1 - t\partial_3 , \\
\tilde P_6 :=
\tilde P_5\frac{t}2p_1 + \twist{P_2}\frac{t-1}3p_1 &= \frac{(4-t^2)t}6p_1^2 + \frac{t^2+t-2}6p_2^2 + \frac{t-1}3p_2 + \frac{t}3p_1^2\partial_2 \\
    &\qquad\qquad\qquad {} + \frac{t-4}6 - \frac{t}2p_1\partial_1 + \frac{2(t-1)}3p_2\partial_2 - \frac{t^2}2p_1\partial_3 , \\
\tilde P_7 :=
\tilde P_6 + \twist{P_2}\frac{t^2-4}3 &= \frac{t^2+t-2}6p_2^2 + \frac{t^3+2t^2-2t-10}6p_2 + \frac{t}3p_1^2\partial_2 \\
    &\qquad\qquad\qquad {} + \frac{2t^2+t-4}6 - \frac{t}2p_1\partial_1 + \frac{2(t-1)}3p_2\partial_2 - \frac{t^2}2p_1\partial_3 + \frac{2(t^2-4)}3\partial_2 .
\end{align*}
Observe how at each line, one can determine precisely
the action of the operator on a monomial~$p_1^{\alpha_1}p_2^{\alpha_2}$
and thus predict the leading monomial of the result
for the monomial ordering refining total degree by~$p_1 > p_2$:
\begin{align*}
\tilde P_5\cdot p_1^{\alpha_1}p_2^{\alpha_2} &= \frac{1-t}3p_1^{\alpha_1+1}p_2^{\alpha_2+1}+\dotsb, \\
\tilde P_6\cdot p_1^{\alpha_1}p_2^{\alpha_2} &= \frac{(4-t^2)t}6p_1^{\alpha_1+2}p_2^{\alpha_2}+\dotsb, \\
\tilde P_7\cdot p_1^{\alpha_1}p_2^{\alpha_2} &= \frac{t^2+t-2}6p_1^{\alpha_1}p_2^{\alpha_2+2}+\dotsb.
\end{align*}
Considering in particular~$\tilde P_7$,
one obtains that degrees with respect to~$p_2$ can be reduced
down to at most~$1$.
Note that the $\tilde P_7\cdot p_1^{\alpha_1}p_2^{\alpha_2}$ luckily
do not reintroduce the variables $p_3$ and~$p_4$.
So at this point, any polynomial~$s \in \bQ(t)[p]$
in an expression $\scalp{F}{sG}$
can be replaced with a linear combination of $1$, $p_1$, $p_2$, and~$p_1p_2$
over~$\bQ(t)$, that is,
with some polynomial confined to a 4-dimensional vector space.
Finally, because $\tilde P_5\cdot 1 = \frac{1-t}3p_1p_2 + \frac{4-t^2}3p_1$,
the monomial~$p_1p_2$ can be replaced with~$p_1$ in such linear combinations,
bringing the finite dimension down to~3.
In the end, for any~$s \in \bQ(t)[p]$,
a sequence of transformations results:
first in an element
$\check s \in \bQ(t) + \bQ(t)p_1 + \bQ(t)p_2$
and elements $\tilde s_j \in \bQ(t)[p]$ for~$j=0,\dots,4$
such that $\scalp{F}{sG} = \scalp{F}{\check sG}$ and
\begin{equation*}
s - \check s = \sum_{i=0}^4 G_i\cdot \tilde s_i
\quad\text{for}\quad
(G_0,\dots,G_4) = (\twist{P_4}, \twist{P_3}, \twist{P_2}, \tilde P_7, \tilde P_5) ;
\end{equation*}
next, because $\tilde P_5$ and~$\tilde P_7$ are in the right ideal,
in elements $\bar s_j \in \bQ(t)[p]$ for~$j=1,\dots,4$ such that
$s - \check s = \twist{P_1}\cdot \bar s_1 + \twist{P_2}\cdot \bar s_2 + \twist{P_3}\cdot \bar s_3 + \twist{P_4}\cdot \bar s_4$.

Eliminating variables one after the other in this presentation
was chosen for the sake of the informal explanation.
In the next section and in our implementation,
we use an optimized elimination strategy that bases more strongly on total degree.

\subsection{Recombining normal forms for a differential equation}

We now explain how the reduction step of the previous section
can be used to derive a differential equation with respect to~$t$
for~$\scalp{F}{G}$.

For any~$i\in\bN$, the identity
$\partial_t^i\cdot\scalp{F}{G} = \scalp{F}{g^iG}$
follows from the definition~$G=\exp(tg)$.
By the reduction of previous section,
the polynomial~$g^i$ can be replaced
with some element~$\check g_i$
from the 3-dimensional vector space $\bQ(t) + \bQ(t)p_1 + \bQ(t)p_2$.
So, the family $\{\check g_0,\check g_1,\check g_2,\check g_3\}$
is obviously linearly dependent over~$\bQ(t)$,
and a linear relation $q_0 \check g_0 + \dots + q_3 \check g_3 = 0$
with~$q_i \in \bQ(t)$
provides a linear differential relation
$(q_0 + q_1 \partial_t + q_2 \partial_t^2 + q_3 \partial_t^3) \cdot \scalp{F}{G} = 0$.

Performing these calculations on our worked example,
we start with $g^0=1$, so that $\check g_0 = 1$ as $1$~is already reduced.
Next,
reducing~$g$ yields $g = \check g_1 + \sum_{i=0}^4 G_i \cdot \tilde s_i$ with
\begin{gather*}
\check g_1 = -\frac{(t^5+2t^4+2t^2+8t-4)}{4(t^2+t-2)t^2}(p_2+1) \\
\text{and}\quad
(\tilde s_0, \dots, \tilde s_4) =
\left(
  \frac14 ,
  \frac{p_1}3 ,
  \frac{p_1^2}{12t} + \frac{(5t-2)p_2}{12t^2} + \frac{4t^2-1}{6t^2} ,
   - \frac{t^2+4t-2}{2t^2(t^2+t-2)} ,
  0
\right) .
\end{gather*}
At this point, a more heavy calculation yields $g^2 = \check g_2 + \sum_{i=0}^4 G_i \cdot \vardbtilde s_i$ with
\begin{align*}
\check g_2 = &-\frac{t^{12} - 14t^{10} - 20t^9 - 36t^8 - 200t^7 - 356t^6 - 48t^5 + 200t^4 - 336t^3 - 240t^2 + 416t - 96}{16(t^2 + t - 2)^2(t - 1)t^4(t + 2)} \\
&- \frac{(t^{13} + 4t^{12} - 16t^{10} - 10t^9 - 36t^8 - 220t^7 - 348t^6 - 48t^5 + 200t^4 - 336t^3 - 240t^2 + 416t - 96)}{16(t^2 + t - 2)^2(t - 1)t^4(t + 2)}p_2
\end{align*}
and quotients~$\vardbtilde s_i$ that we refrain from displaying.
After finding a linear dependency between the~$\check g_i$ over~$\bQ(t)$,
we obtain the annihilating operator
\begin{multline*}
16t^2(t+2)^2(t-1)^2(t^5+2t^4+2t^2+8t-4)\partial_t^2 \\
+(-4t^{13}-16t^{12}+64t^{10}+40t^9+144t^8+880t^7+1392t^6 \\
  \qquad\qquad\qquad+192t^5-800t^4+1344t^3+960t^2-1664t+384)\partial_t \\
-t^4(t^5+2t^4+2t^2+8t-4)^2 .
\end{multline*}
Getting an order~2 less than the dimension~3 could not be predicted.

For efficiency,
the remainders~$\check g_i$ can be obtained in a more incremental way:
the formula
\begin{equation*}
\partial_t^{i+1}\cdot\scalp{F}{G} = \partial_t\cdot\scalp{F}{\check g_iG}
  = \scalp{F}{\partial_t\cdot(\check g_iG)}
  = \scalp{F}{(\check g_i \times g + \partial_t\cdot\check g_i)G}
\end{equation*}
suggests one can obtain~$\check g_{i+1}$
by reducing~$\check g_i \times g + \partial_t\cdot\check g_i$,
which is much smaller than~$g^{i+1}$.
This makes calculations generally faster,
although in the present example
$\check g_1 \times g + \partial_t\cdot\check g_1$
is messier than~$g^2$.

\section{Applicability to various models of graphs}
\label{sec:models}

As we remarked in the introduction, there are many enumeration problems
that can be expressed using the scalar product,
and have the potential to be solved with our strategy.
The computational limits are directly related to the maximal index~$i$ of all $p_i$
that appear in the expressions, and this leaves substantial flexibility.
Although in the work above (namely Section~1.2 and Section~2)
we have focused on the case of simple, loopless graphs,
with only minor modifications of $\serF$ in Eq.~\eqref{eq:Rk-as-scalproduct-of-G}
we can consider graphs with multiple edges, or loops, or both,
as we will prove in Lemma~\ref{lem:scalp-of-f-and-g}.
The form is still an exponential of a polynomial in the~$p_i$.
Similarly, it is straightforward to consider graph classes
where the possible vertex degrees come from a finite set~$K$.
To this end,
it suffices to replace $\exp(t h_k)$ with~$\exp(t(\sum_{j\in K} h_j))$
(as per the lemma again)
and to express the~$h_j$ in the power-sum basis.
For the lemma and future discussions,
we label generalized regular graph models according to three parameters:
\begin{itemize}
\item $e$~encodes the model of allowed edges:
  `se'~is used for graphs with single edges;
  `me'~is used for generalized structures with multiple edges allowed
  (usually called “multigraphs”).
\item $l$~encodes how loops are allowed and counted:
  \begin{itemize}
  \item `ll'~is used for loopless structures, like “graphs” in the usual terminology;
  \item `la'~is used for structures with loops allowed and contributing~$2$ each
  to the degree of a vertex,
  in other words, those models enumerate structures according
  to the number of adjacent half-edges;
  \item `lh'~is used for structures with loops allowed and contributing~$1$ each
  to the degree of a vertex,
  in other words, those models enumerate structures according
  to the number of adjacent edges.
  \end{itemize}
\item $K$~denotes the set of allowed degrees of vertices,
  whether it be counting adjacent edges with `lh'~models
  or counting adjacent half-edges with `la'~models;
  usual $k$-regular graphs are obtained by setting~$K$ to the singleton~$\{k\}$;
  models with~$K$ of larger cardinality allow different vertices of a graph
  to have different degrees as long as they are in~$K$;
  for example, $K = \{1,2,\dots,k\}$ can be used to describe a class of graphs
  with vertex degree bounded by~$k$;
  unless otherwise clear by the context, we make $k = \max K$.
\end{itemize}

\subsection{Theoretical flexibility}

\begin{lemma}\label{lem:scalp-of-f-and-g}
The exponential generating function of a graph model given by some tuple~$(e, l, K)$ is
the scalar product $\scalp{\exp(f)}{\exp(tg)}$
for the polynomials $f$ and~$g$ in $p_1,\dots,p_k$ ($k=\max K$)
defined by Eqs.~\eqref{eq:def-f} and~\eqref{eq:def-g}.
\end{lemma}
  The following proof generalizes
  Eq.~\eqref{eq:Rk-as-scalproduct-of-G} to handle more graph classes,
  and degree restrictions. For each graph
  class we define a symmetric function encoding of
  graphs without degree restrictions, which we shall denote by~$\serF$,
  expressed in the power sum basis (see Table~\ref{tab:sfs}). In all
  cases $\serF$~can be written as an exponential of an infinite sum~$\serf$ of terms in the~$p_i$,
  and the wanted generating function takes the form of a scalar product~$\scalp{\serF}{G}$.
  We show that the coefficient extractor~$G$ has the form
  $G=\exp(tg)$ for a symmetric polynomial $g$. Writing $g$ in the power sum
  basis involves only a finite number of $p_i$, and hence
  $\scalp{\serF}{G} =\scalp{F}{G} $, where $F$ is obtained
  from $\serF$
  by setting all but a finite number of the~$p_i$ to~0, and hence is
  of the form $F=\exp(f)$
  where $f$~is obtained from~$\serf$ in the same way.

\begin{proof}
First we consider the extraction operators and the corresponding series~$G$.
For any symmetric function~$S$,
the coefficient of the monomial $m_\lambda$ in~$S$ is $\scalp{S}{h_\lambda}$.
To get the desired form of~$G$ we use the decomposition
$h_\lambda=h_{\lambda_1}h_{\lambda_2}\dotsm$, and the linearity of the
scalar product.

Thus,
in order to count $k$-regular objects, that is, for the case~$K = \{k\}$,
as above with Eq.~\eqref{eq:Rk-as-scalproduct-of-G}
we have to use $\lambda=(k,\dots,k)$, with $n$~equal parts, for each size~$n$.
This yields the extraction formula
\begin{equation*}
\sum_{n\geq 0} \Scalp{S}{h_{k^n}} \frac{t^n}{n!}=\sum_{n\geq 0} \Scalp{S}{h_k^n \frac{t^n}{n!}} = \scalp{S}{\exp(h_kt)} ,
\end{equation*}
and so~$G = \exp(th_k)$.

In the case where the degree can be from a finite set~$K$ of integers,
we need to extract the coefficients~$\scalp{S}{m_\lambda}$
for partitions~$\lambda$ with all parts in~$K$.
We note that, by the classic correspondence between the exponential function and labelled set constructions,
$G = \exp(t\sum_{k\in K}h_k)$ gives the correct set of monomials with the correct weighting.
Using the change of basis formula $h_n=\sum_{\lambda\vdash n}
\frac{p_\lambda}{z_\lambda}$, we see that written in the power sum
basis, $G$ uses a finite number of $p_i$, the maximum index of which
is at most the maximum element of $K$.

Next, let us consider the symmetric functions~$\serF$ encoding the different unconstrained graph classes
corresponding to each choice for~$(e,l)$.
We can build up the generating function for all six combinations,
and use some basic symmetric-function identities
to express them using power sums.
The results are derived from \cite[Proposition~7.7.4]{stanley_enumerative_1999},
which is proved in a manner similar to Eq.~\ref{eqn:Gequiv} and states
\begin{equation}\label{eq:stanley}
  \prod_{i,j}\frac{1}{1-x_iy_j} =\exp\Biggl(\sum_{n\geq 1} \frac{1}{n} p_n(x)p_n(y)\Biggr)
  \quad \text{and}\quad
  \prod_{i,j} (1+x_iy_j) =\exp\Biggl(\sum_{n\geq 1} \frac{(-1)^{n-1}}{n} p_n(x)p_n(y)\Biggr) ,
\end{equation}
where $p_n(x)$~is the same series $p_n = \sum_i x_i^n$ as before
and $p_n(y)$~is its analogue $\sum_i y_i^n$.
We exploit these equations using two key evaluations.
Setting $y_i=x_i$ in~Eq.~\eqref{eq:stanley} (and hence writing $p_n$ for $p_n(x)$) we get
\begin{equation}\label{eq:prod-xi-xj}
  \prod_{i,j}\frac{1}{1-x_ix_j} =\exp\Biggl(\sum_{n\geq 1} \frac{p_n^2}{n}\Biggr)
  \quad \text{and}\quad
  \prod_{i,j} (1+x_ix_j) = \exp\Biggl(\sum_{n\geq 1} (-1)^{n-1} \frac{p_n^2}{n}\Biggr) .
\end{equation}
Next, setting $y_1=1$ and $y_k=0$, for $k>1$, into Eq.~\eqref{eq:stanley}, thus forcing~$p_n(y) = 1$, for~$n\geq1$, we have
\begin{equation*}
  \prod_{i}\frac{1}{1-x_i}=\exp\Biggl(\sum_{n\geq 1}\frac{p_n}{n}\Biggr)
  \quad \text{and}\quad
  \prod_{i}(1+x_i)=\exp\Biggl(\sum_{n\geq 1}(-1)^{n-1}\frac{p_n}{n}\Biggr) .
\end{equation*}
Remark that for $p_n=p_n(x_1, x_2, \dots)$, we have  $p_n(x^2_1, x^2_2,\dots)=p_{2n}(x_1, x_2,\dots)=p_{2n}$, hence
\begin{equation}\label{eq:prod-xi^2}
  \prod_{i}\frac{1}{1-x^2_i}=\exp\Biggl(\sum_{n\geq 1}\frac{p_{2n}}{n}\Biggr)
  \quad \text{and}\quad
  \prod_{i}(1+x^2_i) =\exp\Biggl(\sum_{n\geq 1} (-1)^{n-1}\frac{p_{2n}}{n}\Biggr) .
\end{equation}
  All six cases can be derived using  these in various products, and
  the results are summarized in Table~\ref{tab:sfs}.
  For example, multiplying Eqs. \eqref{eq:prod-xi-xj} and~\eqref{eq:prod-xi^2} yields
  the squares of the $x$-expressions for (`se',`la') and~(`me',`la'),
  and correspondingly the doubles of the power-sum expressions~$\serf$.
  These calculations give the
  values in Eq.~\eqref{eq:def-f} once we recall that the maximum
  index of a power-sum symmetric function in~$g$ is bounded, thus the scalar
  product will be unchanged if the power sums with indices higher than
  that bound are set to 0.
  Each truncated expression comes in two sums to accommodate the parts in $p_n$ and~$p_{2n}$ separately.

\begin{table}
\begin{tabular}{lll}
Graph type   & $x$-expression~$\serF$ & $\serf$~in power-sum basis\\\hline
\rule{0pt}{0pt} \\[-10pt]
(`se', `ll') & $\prod_{i<j} (1+x_ix_j)$
             & $\sum_{n\geq 1} (-1)^{n-1}\frac{p_n^2-p_{2n}}{2n}$ \\
(`se', `la') & $\prod_{i\leq j} (1+x_ix_j)$
             & $\sum_{n\geq 1} (-1)^{n-1}\frac{p_n^2+p_{2n}}{2n}$ \\
(`se', `lh') & $\prod_{i<j} (1+x_ix_j) \times \prod_{i}(1+x_i)$
             & $\sum_{n\geq 1} (-1)^{n-1}\left(\frac{p_n^2-p_{2n}}{2n} + \frac{p_n}{n}\right)$ \\
(`me', `ll') & $\prod_{i<j}(1-x_ix_j)^{-1}$
             & $\sum_{n\geq 1} \frac{p_n^2-p_{2n}}{2n}$ \\
(`me', `la') & $\prod_{i\leq j}(1-x_ix_j)^{-1}$
             & $\sum_{n\geq 1} \frac{p_n^2+p_{2n}}{2n}$ \\
(`me', `lh') & $\prod_{i<j}(1-x_ix_j)^{-1} \times \prod_{i}(1-x_i)^{-1}$
             & $\sum_{n\geq 1} \left(\frac{p_n^2-p_{2n}}{2n} + \frac{p_n}{n}\right)$
\end{tabular}
\bigskip\bigskip
\caption{Product expressions to encode labelled graphs of the six
  types considered in  Eqs.~\eqref{eq:def-f}. The product
  expression  $\serF$ is equal to the exponential $\exp(\serf)$ where $\serf$~is the summation in the final
  column. The polynomial~$f$ in the hypotheses is equal to
  $\serf$ where all~$p_n$ are set to
  zero, for $n>k$. }
\label{tab:sfs}
\end{table}

\end{proof}

It is worth it to recall that given two combinatorial classes, and differential equations satisfied
by the generating function of each class, we can determine the
differential equations satisfied by both the sum and the product of the two
generating functions.  This sum and product are respectively the
generating functions of the union and the cartesian product of the two
classes.

\subsection{Practical calculations}

We implemented our method as summarized in Algorithm~\ref{algo:method}
and ran it successfully in Maple. Table~\ref{tab:results} presents the results.
All of our calculations are for sets~$K$ included in~$\{1,2,3,4,5,6,7\}$,
allowing several edge and loop variations.
For example,
we have computed the differential equation satisfied
by the set of labelled graphs with degree bounded by $k = 7$,
that is, for $K = \{1,2,3,4,5,6,7\}$,
and the differential equation satisfied
by the set of labelled graphs with degree exactly $k = 7$,
that is, for $K = \{7\}$.

Following up the remark at the end of the previous section,
from our existing results we could for example easily determine the differential
equations satisfied by the set of graphs that are either 5- or
6-regular.  (In contrast to the set of graphs whose vertices are of
degree either 5 or 6, which we can determine directly).

\section{Description of the approach}
\label{sec:approach}

We now provide more formal details on our method,
which will lead to Algorithm~\ref{algo:method}.
Fix a number~$k$ and, again,
write $p = (p_1,\dots,p_k)$ and~$\partial = (\partial_1,\dots,\partial_k)$
as a shorthand.
The number~$k$ is the level of regularity of graphs,
that is, with the $k$ variables in~$p$
we will be able to express the enumerative series of $k$-regular graphs
and variants with regularity bounded by~$k$.

Introduce the Weyl algebra
\begin{equation*}
W_p :=
\bQ\langle p_1,\dots,p_k,\partial_1,\dots,\partial_k;
\partial_i p_j = p_j \partial_i + \delta_{i,j} , \ 1\leq i,j\leq k\rangle ,
\end{equation*}
where $\delta_{i,j}$~is one if and only if~$i=j$, zero otherwise.
Each $\partial_i$~acts on~$\bQ[t][[p]]$
as the usual derivation operator with respect to~$p_i$.
The following relations are easily derived
for any two series $U$ and~$V$ in $\bQ[t][[p]]$ and any~$i \in \{1,\dots,k\}$:
\begin{equation*}
\scalp{p_i U}{V} = \scalp{U}{i\partial_i\cdot V} ,
\qquad
\scalp{i\partial_i\cdot U}{V} = \scalp{U}{p_i V} .
\end{equation*}
By bilinearity and symmetry,
proving these relations reduces indeed to proving the identity
\begin{equation}\label{eq:adj-on-p}
\scalp{p_ip_\lambda}{p_\nu} = \scalp{p_\lambda}{i\partial_i\cdot p_\nu}
\end{equation}
for any $i$, $\lambda$, and~$\nu$.
We prove it for completeness.
First, the identity holds if $\nu$~does not involve~$i$,
both sides being zero.
So we continue assuming $i$~appears in~$\nu$.
Define $\lambda^+$ as the partition obtained
by adjoining $i$ to~$\lambda$
and consider the integers~$r_i$ as in Eq.~\eqref{eq:z-lambda},
so that the analog of Eq.~\eqref{eq:z-lambda} for~$\lambda^+$
is obtained by incrementing~$r_i$.
Therefore,
$z_{\lambda^+}$~is equal to~$z_\lambda r_i (i+1)$.
Define as well $\nu^-$ as the partition obtained by removing $i$ from~$\nu$,
so that $\partial_i\cdot p_\nu = s_i p_{\nu^-}$
where $s_i$~denotes the number of occurrences of~$i$ in~$\nu$.
In particular, $s_i = r_i+1$ if $\nu = \lambda^+$.
Next,
\begin{equation*}
\scalp{p_ip_\lambda}{p_\nu} = \scalp{p_{\lambda^+}}{p_\nu}
  = \delta_{\lambda^+,\nu} z_\lambda i (r_i+1)
  = i (r_i+1) \delta_{\lambda,\nu^-} z_\lambda
  = i (r_i+1) \scalp{p_\lambda}{p_{\nu^-}}
  = i \scalp{p_\lambda}{\partial_i\cdot p_\nu}
\end{equation*}
and Eq.~\eqref{eq:adj-on-p} is proved.
More generally, for any linear differential operator~$L$,
$\scalp{L\cdot U}{V} = \scalp{U}{\adj{L}\cdot V}$,
where the adjoint $\adj{L}$ of~$L$ is the result of applying
the algebra anti-automorphism of~$W_p$
defined by $\adj{p_i} = i\partial_i$ and~$\adj{\partial_i} = i^{-1}p_i$.
This adjoint operation is an involution.
Note that $W_p[t]$~acts on~$\bQ[t][[p]]$ as well,
but we will restrict the use of this action
to right-hand arguments of scalar products.

With Theorem~\ref{thm:proc-is-correct}, we will
only able to prove the correctness of Algorithm~\ref{algo:method},
but not its completeness.
This is why,
we now proceed to progressively develop sufficient properties satisfied by intermediate Gröbner bases in calculations
for Algorithm~\ref{algo:method} to be successful and return correct outputs.

Let $f$ and~$g$ be two polynomials in~$\bQ[p] \setminus \bQ$,
which we do not want to fix at this point to the polynomial needed for the $k$-regular models provided in Lemma~\ref{lem:scalp-of-f-and-g}.
Introduce
\begin{equation*}
 F := \exp(f) \in \bQ[[p]] , \qquad
 G := \exp(tg) \in \bQ[p][[t]] \cap \bQ[t][[p]] , \qquad
 S := \scalp{F}{G} \in \bQ[[t]] .
\end{equation*}
We will write $f_i$ for~$\partial_i\cdot f$
and $g_i$ for~$\partial_i\cdot g$.
Given an element $P \in W_p$, we define its twist $\twist{P} \in W_p[t]$ by
\begin{equation}\label{eq:def-twist}
\twist{P}(p_1,\dots,p_k,\partial_1,\dots,\partial_k) =
\adj{P}(p_1,\dots,p_k,\partial_1+tg_1,\dots,\partial_k+tg_k) .
\end{equation}
Finally, let $H \subset \bQ(t)[p]$ denote the vector space
\begin{equation}\label{eq:vs-im-of-twist}
H := \sum_{P \in \ann(F)} \twist{P}\cdot \bQ(t)[p]
= \sum_{P \in \twist{\ann(F)}} P\cdot \bQ(t)[p] .
\end{equation}

\begin{lemma}\label{lem:H-cancels-scalp}
For any polynomial $h \in H \subset \bQ(t)[p]$,
the scalar product~$\scalp{F}{hG}$ is zero.
\end{lemma}

\begin{proof}
If $P \in W_p$~annihilates~$F$,
then for any $s \in \bQ[p]$,
\begin{equation}\label{eq:move-adj}
0 = \scalp{P\cdot F}{sG} = \scalp{F}{\adj{P}\cdot(sG)} = \scalp{F}{(\twist{P}\cdot s)G} ,
\end{equation}
where $\twist{P}(p_1,\dots,p_k,\partial_1,\dots,\partial_k) \in W_p[t]$ is defined by Eq.~\eqref{eq:def-twist}.
Note that when $P$~runs over the left ideal of annihilating operators of~$F$,
denoted~$\ann(F)$,
the transform~$\adj{P}$ runs over the right ideal $\adj{\ann(F)}$,
and likewise $\twist{P}$ runs over the right ideal $\twist{\ann(F)}$.
The result follows by linearity over~$\bQ(t)$.
\end{proof}

Introduce the derivation operator~$\partial_t$ with respect to~$t$
as well as the Weyl algebra
\begin{equation*}
W_t := \bQ\langle t,\partial_t; \partial_t t = t \partial_t + 1\rangle .
\end{equation*}
Observe
\begin{equation}\label{eq:dt^j-S}
\partial_t^j\cdot S = \scalp{F}{\partial_t^j\cdot G} = \scalp{F}{g^jG} ,
\end{equation}
so that,
as a consequence of Lemma~\ref{lem:H-cancels-scalp},
any annihilator~$Q = \sum_{j=0}^r q_j(t)\partial_t^j \in W_t$ of~$S$
satisfies
\begin{equation*}
0 = Q\cdot S = \scalp{F}{Q\cdot G} = \Scalp{F}{\sum_{j=0}^r q_j g^j G}
= \Scalp{F}{\sum_{j=0}^r q_j (g^j + \ell_j) G}
\end{equation*}
for any polynomials~$\ell_j$ that are
elements of the vector space~$H$ defined by Eq.~\eqref{eq:vs-im-of-twist}.
In what follows, for each~$g^j$ we (implicitly) obtain~$\ell_j$
in such a way that the computed~$g^j + \ell_j$ is “reduced”
and confined in a finite-dimensional $\bQ(t)$-vector space.
This makes it possible to derive the~$q_j$.
The procedure is therefore to deal with Eq.~\eqref{eq:dt^j-S} for each~$j$ separately,
by reducing the coefficient~$g^j$ modulo~$H$.

Although, as we will see, we will only be able to reduce by a subspace of~$H$,
we continue our analysis by expressing the space~$H$ as a finite sum of spaces.
\begin{lemma}\label{lem:H-as-a-finite-sum}
The vector space~$H$ is the sum $\sum_{i=1}^\ell L_i\cdot \bQ(t)[p]$
for any finite family $\{L_i\}_{i=1}^\ell$ generating~$\twist{\ann(F)}$.
\end{lemma}

\begin{proof}
The finite sum $\tilde H := \sum_{i=1}^\ell L_i\cdot \bQ(t)[p]$ is a subspace of~$H$.
Writing any~$P$ of~$\twist{\ann(F)}$ in the form $P = \sum_{i=1}^\ell L_i U_i$
yields the inclusion of the space~$H$ into~$\tilde H$,
and thus the equality~$\tilde H = H$.
\end{proof}

Eq.~\eqref{eq:move-adj}~holds in particular for $P = P_i := i(\partial_i - f_i)$,
in which case $\twist{P}$~is given as~$\twist{P_i}$ in Eq.~\eqref{eq:twist-ann-f}.

To define the reduction that was announced, we proceed by exchanging
the generating family~$\{\twist{P_i}\}$ of~$\twist{\ann(F)}$
for a family $\{G_i\}_{i=1}^\ell$ satisfying the property
that any term~$c m$ to be reduced ($c$~a coefficient, $m$~a monomial)
will be obtained for some $(j, s) \in \{1,\dots,\ell\} \times \bQ(t)[p]$
as the leading monomial of~$G_j\cdot s$,
where leading monomials are decided by some monomial ordering of~$\bQ(t)[p]$
that is compatible with the choice of the family~$\{G_i\}_{i=1}^\ell$.
To make this possible, we ensure that $G_j = m_j + \dotsb \in W_p(t)$
for a monomial~$m_j$ in~$p$,
with the property that, for any $\tilde s \in \bQ(t)[p]$,
the leading monomial of~$m_j \tilde s$ is larger than
the leading monomial of~$(G_j-m_j)\cdot \tilde s$.
In practice, the polynomial~$s$ used to reduce~$c m$ will be set
to the term~$c m/m_j$,
so that $m$~is reduced into~$m - G_j\cdot(c m/m_j)$.
Observing that only finitely many monomials are divisible by none of the~$m_j$
will then ensure the wanted confinement in finite dimension.

The following lemma describes a situation
in which a skew polynomial~$G$ has the property we require from the~$G_j$.
For future reference, we give a name to this property.
\begin{defn}
A skew polynomial~$G \in W_p(t)$ is said to be \emph{dominant}
if it is of the form $G = m + R$ for a non-zero monomial $m \in \bQ[p]$
and some rest~$R$ involving only monomials~$p^\alpha \partial_p^\beta$ for which
$\sum_{i=1}^k (\alpha_i - \beta_i)$~is less than the total degree of~$m$.
\end{defn}

\begin{lemma}\label{lem:lm-is-weighted-highest}
Let $G \in W_p(t)$ be dominant.
Fix a monomial ordering on~$\bQ(t)[p]$ that is graded by total degree.
Then, for any non-zero polynomial~$s$, the leading monomial of~$G\cdot s$
is the product of~$m$ with the leading monomial of~$s$.
\end{lemma}

\begin{proof}
The quantity $\sum_{i=1}^k (\alpha_i - \beta_i)$~is what is added to the degree of any monomial~$q$
when applying~$p^\alpha \partial_p^\beta$ to it.
So because $G$~is dominant and the ordering is graded by total degree,
the leading monomial of~$G\cdot s$ is the product of~$mq$.
The result then follows from the general properties of monomial orderings.
\end{proof}

In general,
we do not know how to ensure the existence of a family $\{G_j\}_{j=1}^\ell$
consisting of dominant~$G_j$,
but for $k$-regular graphs,
a module Gröbner basis calculation will in practice compute such a family for~$k$ up to~$7$,
as we will show in Proposition~\ref{prop:all-Gi-have-property}.
To describe how to engineer the construction of such a family,
we now introduce a few compatible monomial orderings for the structures that we will use.

\begin{defn}\label{def:chosen-orderings}
\begin{enumerate}

\item Let $\prec$ denote a monomial ordering on~$W_p(t)$
  that compares the~$p_i$ by a total degree order
  and that eliminates~$p$
  by making $p_i$~lexicographically higher than~$\partial_j$ for all $i$ and~$j$,
  Here, lexicographically higher means $p_i \succ \partial_j^n$ for all $i$, $j$, and~$n$.

\item Let $\prec_h$ denote the monomial ordering on the free right module $\eta_0 W_p(t) \oplus \eta_1 W_p(t)$ that eliminates~$\eta_1$
  by making $\eta_1$~lexicographically higher than~$\eta_0$ and by sorting monomials in the same~$\eta_i$ by means of~$\prec$.
  Here, $\eta_0$ and~$\eta_1$ are new names denoting elements of a basis and
  lexicographically higher means $p^\alpha \partial^\beta \eta_1 \succ_h p^\gamma \partial^\delta \eta_0$ for all $\alpha$, $\beta$, $\gamma$, and~$\delta$.

\item Let $\prec_p$ denote the monomials ordering on~$\bQ(t)[p]$ that is induced by~$\prec$.
  This is just notation to stress the polynomial situation, as we could use the symbol~$\prec$ as well.

\end{enumerate}
\end{defn}

For each~$i$, write $\twist{P_i} = Q_i(p) + R_i(p,\partial)$,
where $Q_i$~does not involve any~$\partial_j$
and each monomial of~$R_i$ involves at least one~$\partial_j$.
Then, consider $M_i := \eta_1 Q_i + \eta_0 R_i$.
Consider a Gröbner basis for the right\footnote{%
As Maple only computes Gröbner bases for left structures,
the actual computer calculation
computes a Gröbner basis for the left module
generated by the $\adj{Q_i} \eta_1 + \adj{R_i} \eta_0$,
then returns the adjoints $\adj{(Q + R)}$
obtained from the elements $Q \eta_1 + R \eta_0$ of the Gröbner basis
satisfying~$Q \neq 0$.}
module over~$W_p(t)$ generated by the~$M_i$
with respect to~$\prec_\eta$.
Those elements $\eta_1 Q + \eta_0 R$ of the Gröbner basis satisfying~$Q \neq 0$
need not make~$Q+R$ dominant by general properties of Gröbner bases,
but as we start from the dominant elements~$\twist{P_i}$,
the Gröbner basis elements can be hoped to make the~$Q+R$ dominant,
at least if the Gröbner basis calculation does not modify too much the higher monomials of the input polynomials.
Indeed, we are in the nice situation that
the~$Q+R$ are the~$G_i$ we are looking~for
when~$k \leq 6$ for all variant models described in Section~\ref{sec:models}
as well as for the case of regular simple graphs when~$k=7$.
Because the goal of the calculation is to reveal
a zero-dimensional ideal in~$\bQ(t)[p]$,
the module structure over~$W_p(t)$ can be replaced
with a module structure over~$\bQ(t)[p]$,
that is, one would like to consider
only polynomial recombinations guided by the coefficients of~$\eta_1$,
without continuing with non-commutative recombinations
of the coefficients of~$\eta_0$
between generators with zero coefficient with respect to~$\eta_1$.
This is achieved by viewing the~$M_i$ as elements of
a free $\bQ(t)[p]$-module
with a finite basis consisting of $\eta_1$ and some~$\partial^\alpha\eta_0$,
for the same ordering~$\prec_h$.
For all models considered when $2 \leq k \leq 6$,
this has the nice consequence of speeding up the calculation.
For the single model we considered with~$k=7$
we therefore directy used this improvement.

\begin{prop}\label{prop:first-part-correct}
The operators~$G_i$ obtained at Step~\ref{step:optimistic-GB}\@ of Algorithm~\ref{algo:method} are such that
$\sum_{i=1}^\rho G_i W_p(t)$~is a subideal of~$\twist{\ann(F)}$, including the possible degenerate case~0 obtained if~$\rho = 0$.
\end{prop}

\begin{proof}
Given a graph model~$(e, l, K)$, Step~\ref{step:f-and-g}\@ implements the formula announced by Lemma~\ref{lem:scalp-of-f-and-g},
so that $\scalp{\exp(f)}{\exp(tg)}$ is the exponential generating function of the model and $F$~is fixed to~$\exp(f)$.
Next, Step~\ref{step:get-gens}\@ computes generators of the right $W_p(t)$-ideal~$\twist{\ann(F)}$
because the map~$\sharp$ is a linear anti-homomorphism from~$W_p$ to~$W_p(t)$ and the~$P_i$ generate the left $W_p$-ideal~$\ann(F)$.
Steps \ref{step:eta-intro}\@ to~\ref{step:optimistic-GB}\@ produce operators~$G_i$ that are elements of~$\twist{\ann(F)}$.
This proves the result.
\end{proof}

\begin{algorithm}
\begin{quote}
\rule{.85\textwidth}{.1em}

\noindent
{\bfseries Input:}
  a graph model $(e,l,K)$, \\
  \phantom{mmmm}where
  $e \in \{\text{`se'},\text{`me'}\}$,
  $l \in \{\text{`ll'},\text{`la'},\text{`lh'}\}$,
  $K \subset \bN_{>0}$.

\noindent
{\bfseries Output:}
  an operator of minimal order in~$\partial_t$ \\
  \phantom{mmmm}that cancels the exponential generating function of the model.\\[0mm]\noindent
\rule{.85\textwidth}{.1em}
\begin{enumerate}[a.]
\item\label{step:f-and-g}
   Set $k = \max K$, then compute $f$ and~$g$ by the formulas
  \begin{align}
    \label{eq:def-f}
    f &= \sum_{i=1}^k
            \bigl((-1)^{i+1} \delta_{e,\text{`se'}} + \delta_{e,\text{`me'}}\bigr)
            \left( \frac{p_i^2}{2i} + \delta_{l,\text{`lh'}} \frac{p_i}{i} \right) \\
      &\notag \qquad\qquad\qquad
       - (-1)^{\delta_{l,\text{`la'}}} \sum_{i=1}^{\lfloor k/2\rfloor}
            \bigl((-1)^{i+1} \delta_{e,\text{`se'}} + \delta_{e,\text{`me'}}\bigr)
            \frac{p_{2i}}{2i} , \\
    \label{eq:def-g}
    g &= \sum_{j \in K} \sum_{\lambda\vdash j} \frac{p_\lambda}{z_\lambda} ,
  \end{align}
  where the expression~$\delta_{x,y}$ is equal to~1 if~$x=y$ and to~0 otherwise.
\item\label{step:get-gens}
  Get generators of the right $W_p(t)$-ideal~$\twist{\ann(F)}$
  by computing $\twist{P_i}$ for $1\leq i\leq k$ by the formula
\begin{equation}\label{eq:twist-ann-f}
\twist{P_i} =
\twist{(i\partial_i - if_i)} = p_i - if_i(\partial_1+tg_1,2(\partial_2+tg_2),\dots,k(\partial_k+tg_k)) .
\end{equation}
Here, the right-hand side is obtained by the non-commutative substitution
of $p_1$ with~$\partial_1+tg_1$, of $p_2$ with~$2(\partial_2+tg_2)$, \dots, of $p_k$ with~$k(\partial_k+tg_k)$,
in the polynomial~$f_i = f_i(p_1,\dots,p_k)$.
\item\label{step:eta-intro}
  Transform each~$\twist{P_i}$ by the map
  \[ \sum_\alpha c_\alpha \partial^\alpha \mapsto c_0 \eta_1 + \sum_{\alpha\neq0} c_\alpha \partial^\alpha \eta_0 \]
  to get a system of generators of a right $\bQ(t)[p]$-module
  with basis $\{\eta_1\} \cup \{\partial^\alpha \eta_0\}_\alpha$.
  Here,
  $\alpha$~ranges in the finite set of non-zero exponents involved in the~$\twist{P_i}$.
\item\label{step:GB-elim-eta}
  Compute a Gröbner basis of this right module for an ordering~$\prec_\eta$ that makes
  $\eta$ lexicographically higher than~$p$
  and $p$ lexicographically higher than~$\partial$,
  using Definition~\ref{def:chosen-orderings}.
\item\label{step:optimistic-GB}
  Obtain elements $G_1,\dots,G_\rho$ of~$W_p(t)$ by setting $\eta_1 = \eta_0 = 1$
  in those elements of the Gröbner basis
  that involve~$\eta_1$ with a non-zero coefficient,
  then write each~$G_i$ in the form~$Q_i(p) + R_i(p,\partial)$
  where each monomial of~$R_i$ involves at least one~$\partial_j$.
\item\label{step:test-dim}
  If the polynomial ideal~$I = (Q_1,\dots,Q_\rho)$ of~$\bQ(t)[p]$ has positive dimension,
  then return `FAIL',
  else determine
  the monomials $p^{\beta_1},\dots,p^{\beta_\delta}$ under the stair of~$I$ with respect to~$\prec_p$.
\item\label{step:reductions}
  Set $\check g_0 = 1$, then for~$i$ from~2 to~$\delta$,
  set $\check g_i = \red(g \check g_{i-1} + \partial_t\cdot \check g_{i-1}, (G_i)_{i=1}^\rho, \prec)$.
\item\label{step:matrix}
  Compute the matrix~$M$ with rows indexed by~$0\leq i\leq\delta$
  and columns indexed by~$1\leq j\leq\delta$,
  whose entry at position~$(i,j)$
  is the coefficient of~$p^{\beta_j}$ in~$G_i$.
\item\label{step:kernel}
  Compute a basis of the left kernel of~$M$,
  then combine its elements to obtain a non-zero row vector
  $(q_0,\dots,q_\delta) \in \bQ(t)^{\delta+1}$
  with maximal number of zeros to the right.
\item\label{step:return-ODE}
  Return $q_0 + q_1\partial_t + \dots + q_\delta\partial_\delta$.
\end{enumerate}
\rule{.85\textwidth}{.1em}
\end{quote}
\caption{\label{algo:method}
Outline of the method. Uses the reduction $\red({\cdot})$ of {\sc{Algorithm}}~\ref{algo:reduction}.}
\end{algorithm}

\begin{prop}\label{prop:all-Gi-have-property}
For each model in Table~\ref{tab:results}:
\begin{enumerate}
\item\label{it:Q_i+R_i-ok}
  the elements $Q_i(p) + R_i(p,\partial)$ of the Gröbner basis obtained at Step~\ref{step:optimistic-GB}\@ of Algorithm~\ref{algo:method}
  are all dominant and so satisfy the property of~$G$ in Lemma~\ref{lem:lm-is-weighted-highest},
  for~$m=Q_i(p)$.
\item\label{it:I-has-dim-0}
  the ideal~$I$ generated by the~$Q_i(p)$ at Step~\ref{step:test-dim}\@ has dimension zero.
\end{enumerate}
\end{prop}

\begin{proof}
For each model,
the proof is by inspection after computing the Gröbner basis at Step~\ref{step:GB-elim-eta}\@:
testing Point~(\ref{it:Q_i+R_i-ok}) consists in comparing each monomial of~$Q_i+R_i$ with the leading monomial~$m_i$ of~$Q_i$;
testing Point~(\ref{it:I-has-dim-0}) is done by computing the dimension of the commutative polynomial ideal by a classical algorithm,
after observing that the commutative polynomials~$Q_i$ obtained as parts of the~$Q_i+R_i$ are already a Gröbner basis
for the ordering induced on~$\bQ(t)[p]$ by the ordering~$\prec_\eta$ used at Step~\ref{step:GB-elim-eta}\@
\end{proof}

The polynomials~$G_i = Q_i+R_i$ obtained at Step~\ref{step:optimistic-GB}\@ of Algorithm~\ref{algo:method} will be used
to reduce polynomials to a finite-dimensional vector space at Step~\ref{step:reductions}\@
So we first analyse this reduction separately.
To this end,
we stress a possibly confusing fact:
although $\twist{\ann(F)}$~is a right-module,
we make it act to the left of polynomials during reductions.

\begin{algorithm}
\begin{quote}
\rule{.85\textwidth}{.1em}

\noindent
{\bfseries Input:}
  a polynomial~$s \in \bQ(t)[p]$ to be reduced
  and dominant operators $G_1,\dots,G_\rho$ of~$W_p(t)$.

\noindent
{\bfseries Output:}
  a polynomial~$\check s \in \bQ(t)[p]$.\\[0mm]\noindent
\rule{.85\textwidth}{.1em}
\begin{enumerate}[a.]
\item For each~$i$, set $m_i$ to the leading monomial of~$G_i$.
\item If no monomial of~$s$ is divisible by any of the~$m_i$, return~$s$.
\item Set $m$ to the maximal monomial in~$s$ that is divisible by some~$m_i$
  and choose $j$ such that $m_j$~divides~$m$.
\item Set $c$~to the coefficient of~$m$ in~$s$
  and $c_j$~to the leading coefficient of~$G_j$.
\item\label{step:red-do-it}
  Set $t$ to the term~$\frac c{c_j} \frac m{m_j}$
  and return $\red(s - G_j\cdot t, (G_i)_{i=1}^\rho)$.
\end{enumerate}
\rule{.85\textwidth}{.1em}
\end{quote}
\caption{\label{algo:reduction}
Reduction $\red({\cdot})$ used by {\sc{Algorithm}}~\ref{algo:method}.}
\end{algorithm}

\begin{prop}\label{prop:red-correct}
The reduction algorithm, Algorithm~\ref{algo:reduction}, terminates
if the ideal~$I$ generated by the leading monomials with respect to~$\prec$ of the dominant inputs~$G_i$
is zero-dimensional.
It returns a polynomial~$\check s \in \bQ(t)[p]$ whose monomials are under the stairs of~$I$ with respect to~$\prec_p$.
The difference~$s - \check s$ is in~$\sum_{i=1}^\rho G_i\cdot \bQ(t)[p]$,
and therefore it is an element of the vector space~$H$ defined by Eq.~\eqref{eq:vs-im-of-twist}.
\end{prop}

\begin{proof}
The term~$t$ at Step~\ref{step:red-do-it}\@ is chosen so that leading terms satisfy:
\[ \lt(G_j\cdot t) = \lt(\lt(G_j)\cdot t) = \lt(c_j m_j\cdot t) = \lt(c_j m_j t) = \lt(c m) = c m , \]
where the first equality is by the assumed dominant character of~$G_j$ and by Lemma~\ref{lem:lm-is-weighted-highest},
and where the third equality is because $m_j$~does not involve~$\partial$.
So all monomials appearing in the difference~$s - G_j\cdot t$ and susceptible of reduction by the~$m_i$ are less than~$m$.
Because of the 0-dimensionality of~$I$, the recursive calls to~$\red({\cdot})$ terminate.
The output has no monomial reducible by any~$m_i$, hence all its monomials are under the stairs of~$I$.
Finally, each~$G_j\cdot t$ considered by some recursive call is in~$H$, hence so is~$\check s - s$.
\end{proof}

We note that the choice of~$\prec_p$ to be induced by~$\prec_h$
and that $\prec_h$~reduces to a graded order on monomials in~$p$
is not for termination, but for efficiency.

\begin{theorem}\label{thm:proc-is-correct}
Algorithm~\ref{algo:method} is correct, that is,
if it terminates, then either this is by giving up, returning `FAIL' at Step~\ref{step:test-dim}\@,
or this is by outputting at Step~\ref{step:return-ODE}\@ a differential equation
that annihilates the scalar product~$S$.
\end{theorem}

\begin{proof}
By Proposition~\ref{prop:first-part-correct}, we obtain that, after Step~\ref{step:optimistic-GB}\@,
$\sum_{i=1}^\rho G_i W_p(t)$-is a subideal of~$\twist{\ann(F)}$.
So, if the algorithm terminates without returning `FAIL' at Step~\ref{step:test-dim}\@,
then the ideal~$I$, which is generated by the~$Q_i = G_i-R_i$, must have dimension~0.
Proposition~\ref{prop:red-correct} shows that $\red({\cdot})$~modifies its input~$s$
by adding to it an element~$\check s - s$ of the vector space~$H$ defined by Eq.~\eqref{eq:vs-im-of-twist}.
So, Step~\ref{step:reductions}\@ is an incremental calculation of the~$\partial_t^i\cdot S$ defined by Eq.~\eqref{eq:dt^j-S},
as justified by Lemma~\ref{lem:H-cancels-scalp}.
The final steps compute a non-trivial $\bQ(t)$-linear relation between the~$\partial_t^i\cdot S$.
The corresponding differential equation is finally returned.
\end{proof}

We emphasize that we do not claim that
the case of 0-dimensionality of the ideal~$I$ at Step~\ref{step:test-dim}\@ implies that termination of the subsequent steps.
In fact, this is just a necessary condition for those steps to be well defined.
A sufficient condition for termination, knowing that the dimension of~$I$ is~0, is that for each~$i$ in~$\{1,\dots,\rho\}$
the (total) degree of~$Q_i$ in~$p$ is larger than the (partial) degree of~$R_i$ in~$p$
(cf. Lemma~\ref{lem:lm-is-weighted-highest}).
This is a condition that we have observed in all of our experiments.

Because we observe that our implementation of Algorithm~\ref{algo:method}
terminates and outputs a differential equation on the models listed in Table~\ref{tab:results},
we get the following corollary.

\begin{corollary}\label{cor:table-is-correct}
For each model in Table~\ref{tab:results}, there exist:
\begin{enumerate}
\item\label{item:lde-is-correct}
  a linear differential equation annihilating the generating function of the model, with order given in column~`$\partial_t$',
\item a linear recurrence equation annihilating the enumerative sequence of the model, with order given in column~`$\partial_n$ minimized',
  the latter being the minimal order of a recurrence if it is not starred.
\end{enumerate}
\end{corollary}

\begin{proof}
Only the minimality of recurrence order requires a proof:
this is by the correctness of van Hoeij's \texttt{LREtools[MinimalRecurrence]} implementation,
whose proven method is described in his student Zhou's PhD thesis~\cite[Chapter~6]{Zhou-2022-AFL}.
\end{proof}

\section{No computation of initial conditions is needed}

For all classes,
after computing the ODE one readily proves by observation
that it possesses the only exponent~$n=0$ at~$t=0$.
Consequently,
the series solutions form a $1$-dimensional vector space,
for which a possible basis is the family with only entry
our combinatorial series,
e.g., the singleton family~$(R^{(k)}(t))$ in the $k$-regular case.
Note that the empty graph is $k$-regular for any~$k$,
implying the identity $R^{(k)}(t) = 1 + O(t)$,
and an analogue result holds for all degree sets~$K$.
Converting the ODE to a recurrence relation
satisfied by the coefficient sequence~$(c_n)_{n\in\bZ}$
of any of its series solution $\sum_{n\in\bZ} c_n t^n$
and forcing $c_0 = 1$ and $c_n = 0$ for all~$n < 0$
therefore uniquely determines all~$c_n$ for~$n > 0$.
This observation generalizes to
any of the models of edges and loops
presented in Section~\ref{sec:models}.

So, a complete proof of correctness of the method
for computing an ODE satisfied by the scalar product,
together with the observation above, makes it unnecessary to apply a
resource-consuming calculation of first terms of the series.
Nonetheless, for all the classes
we did verify our series solution by direct computation of scalar product.
For example for simple $k$-regular graphs,
we directly determined~$r^{(k)}_n$ for
$n\leq 1284$ when~$k=3$,
$n\leq 216$ when~$k=4$,
$n\leq 90$ when~$k=5$,
$n\leq 46$ when~$k=6$, and
$n\leq 31$ when~$k=7$.
Furthermore, McKay provided values for~$k=5$ and~$n\leq 600$, all
consistent with our computations.

\section{Minimal recurrences}

Table~\ref{tab:results} shows that in many instances
the obtained ODE has low order and high degree,
leading to a recurrence of high order and low degree.
This is one possible motivation for searching for recurrences of lower orders.
To this end, we have used Mark van Hoeij's Maple implementation
of his algorithm to reduce the order of recurrences
satisfied by a specific solution to a given initial recurrence%
\footnote{The documentation of the procedure promises to achieve the minimal order,
but to the best of our knowledge no formal publication is available yet.},
which is available in Maple as \texttt{LREtools[MinimalRecurrence]}.

For roughly half of the models, we could reduce the order.
For a small quarter of the list, minimizing was too much calculation.
Fortunately in those cases,
the generating series is even,
the initial recurrence relates every second term of the sequence,
and a change of indexes to consider the sub-sequence of even terms
led to a a recurrence that Maple could reduce.
These are marked with a star in the table.

However, all models with~$k=\max K=6$ and loops described by~$l={}$`lh',
no recurrence of lower order could be found.

In all relevant cases, the recurrence of reduced order is larger
than the initial recurrence,
as the degree increase outbalances the order decrease.
This is amplified by an increase in the size of integers that occur in the reduced recurrence.
For example, for the model defined by $e={}$`se', $l={}$`ll', $K=\{6\}$,
the product ``degree${}\times{}$order'' is multiplied by~$4.8$,
raising from~$876$ to~$4176$,
while the length of the longest integer raises from 50 to~236.
As a consequence, computing first terms as with Maple's command \texttt{gfun[rectoproc]}
results in slower calculations with the reduced recurrence:
for the same example, computing up to the 1000th term takes more than five times
as much with the reduced recurrence,
with times raising from below 30~seconds to above 150~seconds.

\section{Conclusion}

\subsection{Computational considerations}
For each graph class we considered with degrees bounded by~$6$,
it did not require more than 10~minutes to determine the differential
equation satisfied the generating function.
In contrast, for models with~$k = \max K = 7$,
the very same implementation requires hours to terminate
(between $4.5$ and~$\simeq30$).
For example, for the model defined by $e={}$`se', $l={}$`ll', $K=\{7\}$,
the time breaks down as follows:
a Gröbner basis that supports the reduction by the vector space~$H$ of Eq.~\eqref{eq:vs-im-of-twist}.
can be obtained in less than $1$~hour
(Steps \ref{step:f-and-g}\@ to~\ref{step:test-dim}\@ in Algorithm~\ref{algo:method}).
From this, one can predict
that reduced forms of scalar products will be confined in dimension~20.
Twenty successive reductions are then performed,
taking longer and longer,
for a total duration of about $9$~hours
(Step~\ref{step:reductions}\@).
After this, the linear algebra
(Steps \ref{step:matrix}\@ and~\ref{step:kernel}\@)
is comparatively fast.
The resulting ODE
satisfied by the generating series~$R^{(7)}(t)$ for 7-regular graphs
has order~20.

Of course the natural question to ask is \emph{What about $k>7$?}
In that respect, ongoing discussions with Hadrien Brochet
have led to promising observations that could speed up calculations
and hopefully get $k=8$.

Finally, the generating function for all regular graphs is not
D-finite. Are there properties of the presented ODEs
that can help us understand
if the generating function of all regular graphs is differentially
algebraic or not, and if so how to find the differential equation?

\subsection{Combinatorial considerations}
We have used symmetric function identities to get an expression in the
power-sum basis for the graph generating functions.  In fact, there is
a well developed combinatorial theory to get these expressions
directly using cycle index series, and other machinery from Species
theory. The underlying combinatorial framework is developed and
rigourous in~\cite{mendez_multisets_1993} and is interpreted for this
context in \cite{mishna_regularity_2018}. Roughly, this means that we
can express the generating functions of a wide family of combinatorial
objects such as hypergraphs, and weighted graphs under restrictions of
graph degree in a form similar to those expressed in
Lemma~\ref{lem:scalp-of-f-and-g}. Furthermore, the $f$ are very simply
deduced using the plethysm operator of symmetric functions.

The symmetric function $\prod_{i\geq 1}\frac{1}{1-x_i}\prod_{1\leq i<j}\frac{1}{1-x_ix_j}$ encodes graphs models with $e={}$`me', $l={}$`ll'.
This is a very well studied symmetric function, as it also
encodes semi-standard Young Tableaux by their content
\cite[Corollary~7.13.8]{stanley_enumerative_1999}.
Thus, the method presented in this work determines generating
functions of semi-standard Young tableaux with restrictions on the
number of times each number appears as an entry. For example, $\langle
\prod_{i\geq 1}\frac{1}{1-x_i}\prod_{1\leq i<j}\frac{1}{1-x_ix_j},
h_1^n\rangle$ is the number of standard Young
tableaux on $n$ boxes with entries $1, 2, \dots, n$, each appearing
exactly once. Some of our results then double as recurrences for the number of semi-standard Young tableaux on~$n$ boxes where the number of times an entry appears comes from a finite set~$K$.

\section{Acknowledgements}
The first author's work was supported in part by the French ANR grant \href{https://mathexp.eu/DeRerumNatura/}{\emph{De rerum natura}} (ANR-19-CE40-0018)
and by the French-Austrian ANR-FWF grant EAGLES (ANR-22-CE91-0007 \& FWF-I-6130-N).
The second author's work was supported by Discovery Grant
RGPIN-2017-04157 from the National Science and Engineering Research
Council of Canada. The authors are grateful to Nick Wormald for
comments on this work, and also to Brendan McKay for motivating our interest
in the problem; for providing some counting
sequence data for comparision purposes; and for keen observations that arose
in the study of the output recurrences.

\bibliographystyle{plain}
\bibliography{references}

\begin{thebibliography}{10}

\bibitem{BostanChenChyzakLi-2010-CCT}
Alin Bostan, Shaoshi Chen, Frédéric Chyzak, and Ziming Li.
\newblock Complexity of creative telescoping for bivariate rational functions.
\newblock In {\em ISSAC'10: Proceedings of the 2010 International Symposium on
  Symbolic and Algebraic Computation}, pages 203--210, New York, NY, USA, 2010.
  ACM.

\bibitem{BostanChyzakLairezSalvy-2018-GHR}
Alin Bostan, Frédéric Chyzak, Pierre Lairez, and Bruno Salvy.
\newblock Generalized {H}ermite reduction, creative telescoping and definite
  integration of {D}-finite functions.
\newblock In Éric Schost, editor, {\em ISSAC'18}, pages 95--102. ACM Press,
  2018.

\bibitem{chyzak_effective_2005}
Frédéric Chyzak, Marni Mishna, and Bruno Salvy.
\newblock Effective scalar products of {D}-finite symmetric functions.
\newblock {\em Journal of Combinatorial Theory. Series A}, 112(1):1--43, 2005.

\bibitem{gessel_enumerative_1987}
Ira Gessel.
\newblock Enumerative applications of symmetric functions.
\newblock {\em Séminaire Lotharingien de Combinatoire [electronic only]},
  17:B17a--17, 1987.
\newblock Publisher: Universität Wien, Fakultät für Mathematik.

\bibitem{gessel_symmetric_1990}
Ira~M. Gessel.
\newblock Symmetric functions and {P}-recursiveness.
\newblock {\em Journal of Combinatorial Theory. Series A}, 53(2):257--285,
  1990.

\bibitem{goulden_hammond_1983}
I.~P. Goulden, D.~M. Jackson, and J.~W. Reilly.
\newblock The {Hammond} series of a symmetric function and its application to
  {P}-recursiveness.
\newblock {\em Society for Industrial and Applied Mathematics. Journal on
  Algebraic and Discrete Methods}, 4(2):179--193, 1983.

\bibitem{gropp_enumeration_1992}
Harald Gropp.
\newblock Enumeration of regular graphs 100 years ago.
\newblock In {\em Discrete {Mathematics}}, volume 101, pages 73--85. 1992.
\newblock ISSN: 0012-365X,1872-681X Journal Abbreviation: Discrete Math.

\bibitem{Lipshitz-1989-DFP}
L.~Lipshitz.
\newblock {$D$}-finite power series.
\newblock {\em Journal of Algebra}, 122(2):353--373, 1989.

\bibitem{mendez_multisets_1993}
Miguel M\'endez.
\newblock Multisets and the combinatorics of symmetric functions.
\newblock {\em Adv. Math.}, 102(1):95--125, 1993.

\bibitem{mishna_automatic_2007}
Marni Mishna.
\newblock Automatic enumeration of regular objects.
\newblock {\em Journal of Integer Sequences}, 10(5):Article 07.5.5, 18, 2007.

\bibitem{mishna_regularity_2018}
Marni Mishna.
\newblock Regularity in weighted graphs a symmetric function approach.
\newblock {\em Contributions to Discrete Mathematics}, 13(2):32--44, 2018.

\bibitem{oeis}
{OEIS Foundation Inc.}
\newblock The {O}n-{L}ine {E}ncyclopedia of {I}nteger {S}equences.
\newblock Published electronically at \url{http://oeis.org}.
\newblock Accessed on July 23, 2024.

\bibitem{read_enumeration_1959}
R.~C. Read.
\newblock The enumeration of locally restricted graphs. {I}.
\newblock {\em The Journal of the London Mathematical Society}, 34:417--436,
  1959.

\bibitem{read_number_1980}
R.~C. Read and N.~C. Wormald.
\newblock Number of labeled 4-regular graphs.
\newblock {\em Journal of Graph Theory}, 4(2):203--212, 1980.

\bibitem{read_some_1959}
Ronald~C Read.
\newblock {\em Some enumeration problems in graph theory}.
\newblock PhD thesis, University of London (University College of the West
  Indies), 1959.

\bibitem{stanley_differentiably_1980}
R.~P. Stanley.
\newblock Differentiably finite power series.
\newblock {\em European Journal of Combinatorics}, 1(2):175--188, 1980.

\bibitem{stanley_enumerative_1999}
Richard~P. Stanley.
\newblock {\em Enumerative combinatorics. {V}ol. 2}, volume~62 of {\em
  Cambridge Studies in Advanced Mathematics}.
\newblock Cambridge University Press, Cambridge, 1999.
\newblock With a foreword by Gian-Carlo Rota and appendix 1 by Sergey Fomin.

\bibitem{wormald_asymptotic_2018}
Nicholas Wormald.
\newblock Asymptotic enumeration of graphs with given degree sequence.
\newblock In {\em Proceedings of the {International} {Congress} of
  {Mathematicians}: {Rio} de {Janeiro} 2018}, pages 3245--3264. World
  Scientific, 2018.

\bibitem{Zhou-2022-AFL}
Yi~Zhou.
\newblock {\em Algorithms for factoring linear recurrence operators}.
\newblock PhD thesis, Florida State University, 2022.

\end{thebibliography}


\appendix

\section{Differential equations and recurrences relations}

Table~\ref{tab:results} gathers information related to computations
we performed for a list of models:
\begin{itemize}
\item parameters “edges”, “loops”, and $k$'s are
  as described in Section~\ref{sec:models};
\item the obtained ODE has order provided in column~$\partial_t$ and
  its polynomial coefficients have degrees bounded by the number in column~$t$;
\item a first recurrence on the number of graphs of size~$n$
  is directly obtained by translating the ODE
  by Maple's \texttt{gfun[diffeqtorec]};
  it has order provided in column~$\partial_n$
  and its polynomial coefficients have degrees bounded by the number in column~$n$;
\item in the majority of models, the first recurrence could be minimized
  by \texttt{LREtools[MinimalRecurrence]},
  leading to a new pair of columns $n$ and~$\partial_n$;
  the minimal order is proved unless it is starred,
  meaning that the minimized recurrence is only of minimal order
  among recurrences on even terms;
\item the corresponding calculation is done in the time of column “time”,
  measured in seconds.
\end{itemize}

\begin{longtable}
{rrr|rr|rrrr|r}
KILLED & LINE!!!! \kill
\caption[An optional table caption ...]{%
Models up to~$k=7$.
See description in Section~\ref{sec:models}.
All timings obtained on the same computer (with AMD EPYC 9754 processor).
\label{tab:results}}\\
\hline
edges & loops & $k$'s & $t$ & $\partial_t$ & $n$ & $\partial_n$ & $n$ & $\partial_n$ & time \\
& & & & & \multicolumn{2}{c}{from ODE} & \multicolumn{2}{c|}{minimized} & \\
\hline
\endfirsthead
\caption[]{Models up to $k=7$ (continued)}\\
\hline
edges & loops & $k$'s & $t$ & $\partial_t$ & $n$ & $\partial_n$ & $n$ & $\partial_n$ & time \\
& & & & & \multicolumn{2}{c}{from ODE} & \multicolumn{2}{c|}{minimized} & \\
\hline
\endhead
\hline
& & & & & \multicolumn{2}{c}{from ODE} & \multicolumn{2}{c|}{minimized} & \\
edges & loops & $k$'s & $t$ & $\partial_t$ & $n$ & $\partial_n$ & $n$ & $\partial_n$ & time \\
\hline
\endfoot
\hline
\endlastfoot
se & ll & 2 & 2 & 1 & 1 & 3 & 1 & 3\phantom{${}^*$} & 0.04 \\
me & ll & 2 & 2 & 1 & 1 & 3 & 1 & 3\phantom{${}^*$} & 0.05 \\
se & la & 2 & 2 & 1 & 1 & 3 & 1 & 3\phantom{${}^*$} & 0.05 \\
me & la & 2 & 2 & 1 & 1 & 3 & 1 & 3\phantom{${}^*$} & 0.05 \\
se & lh & 2 & 3 & 1 & 1 & 4 & 1 & 4\phantom{${}^*$} & 0.05 \\
me & lh & 2 & 3 & 1 & 1 & 4 & 1 & 4\phantom{${}^*$} & 0.05 \\
se & ll & 1,2 & 3 & 1 & 1 & 4 & 1 & 4\phantom{${}^*$} & 0.05 \\
me & ll & 1,2 & 3 & 1 & 1 & 4 & 1 & 4\phantom{${}^*$} & 0.05 \\
se & la & 1,2 & 3 & 1 & 1 & 4 & 1 & 4\phantom{${}^*$} & 0.07 \\
me & la & 1,2 & 3 & 1 & 1 & 4 & 1 & 4\phantom{${}^*$} & 0.05 \\
se & lh & 1,2 & 3 & 1 & 1 & 4 & 1 & 4\phantom{${}^*$} & 0.06 \\
me & lh & 1,2 & 3 & 1 & 1 & 4 & 1 & 4\phantom{${}^*$} & 0.06 \\
se & ll & 3 & 11 & 2 & 2 & 12 & 4 & 8${}^*$ & 0.08 \\
me & ll & 3 & 11 & 2 & 2 & 12 & 4 & 8${}^*$ & 0.08 \\
se & la & 3 & 11 & 2 & 2 & 12 & 4 & 8${}^*$ & 0.09 \\
me & la & 3 & 11 & 2 & 2 & 12 & 4 & 8${}^*$ & 0.07 \\
se & lh & 3 & 11 & 2 & 2 & 12 & 6 & 8\phantom{${}^*$} & 0.09 \\
me & lh & 3 & 11 & 2 & 2 & 12 & 6 & 8\phantom{${}^*$} & 0.1 \\
se & ll & 1,3 & 11 & 2 & 2 & 12 & 4 & 8${}^*$ & 0.09 \\
me & ll & 1,3 & 11 & 2 & 2 & 12 & 4 & 8${}^*$ & 0.07 \\
se & la & 1,3 & 11 & 2 & 2 & 12 & 4 & 8${}^*$ & 0.09 \\
me & la & 1,3 & 11 & 2 & 2 & 12 & 3 & 8${}^*$ & 0.07 \\
se & lh & 1,3 & 11 & 2 & 2 & 12 & 6 & 8\phantom{${}^*$} & 0.07 \\
me & lh & 1,3 & 11 & 2 & 2 & 12 & 6 & 8\phantom{${}^*$} & 0.09 \\
se & ll & 2,3 & 11 & 2 & 2 & 12 & 6 & 8\phantom{${}^*$} & 0.08 \\
me & ll & 2,3 & 11 & 2 & 2 & 12 & 6 & 8\phantom{${}^*$} & 0.1 \\
se & la & 2,3 & 11 & 2 & 2 & 12 & 6 & 8\phantom{${}^*$} & 0.08 \\
me & la & 2,3 & 11 & 2 & 2 & 12 & 6 & 8\phantom{${}^*$} & 0.08 \\
se & lh & 2,3 & 11 & 2 & 2 & 12 & 6 & 8\phantom{${}^*$} & 0.08 \\
me & lh & 2,3 & 11 & 2 & 2 & 12 & 6 & 8\phantom{${}^*$} & 0.08 \\
se & ll & 1,2,3 & 11 & 2 & 2 & 12 & 6 & 8\phantom{${}^*$} & 0.07 \\
me & ll & 1,2,3 & 11 & 2 & 2 & 12 & 6 & 8\phantom{${}^*$} & 0.1 \\
se & la & 1,2,3 & 11 & 2 & 2 & 12 & 6 & 8\phantom{${}^*$} & 0.08 \\
me & la & 1,2,3 & 11 & 2 & 2 & 12 & 6 & 8\phantom{${}^*$} & 0.07 \\
se & lh & 1,2,3 & 11 & 2 & 2 & 12 & 6 & 8\phantom{${}^*$} & 0.08 \\
me & lh & 1,2,3 & 11 & 2 & 2 & 12 & 6 & 8\phantom{${}^*$} & 0.08 \\
se & ll & 4 & 14 & 2 & 2 & 15 & 7 & 10\phantom{${}^*$} & 0.2 \\
me & ll & 4 & 14 & 2 & 2 & 15 & 7 & 10\phantom{${}^*$} & 0.2 \\
se & la & 4 & 14 & 2 & 2 & 15 & 7 & 10\phantom{${}^*$} & 0.19 \\
me & la & 4 & 14 & 2 & 2 & 15 & 6 & 10\phantom{${}^*$} & 0.17 \\
se & lh & 4 & 30 & 3 & 3 & 31 & 18 & 16\phantom{${}^*$} & 0.19 \\
me & lh & 4 & 29 & 3 & 3 & 30 & 17 & 16\phantom{${}^*$} & 0.19 \\
se & ll & 1,4 & 29 & 3 & 3 & 30 & 17 & 16\phantom{${}^*$} & 0.29 \\
me & ll & 1,4 & 29 & 3 & 3 & 30 & 17 & 16\phantom{${}^*$} & 0.28 \\
se & la & 1,4 & 29 & 3 & 3 & 30 & 17 & 16\phantom{${}^*$} & 0.2 \\
me & la & 1,4 & 29 & 3 & 3 & 30 & 17 & 16\phantom{${}^*$} & 0.18 \\
se & lh & 1,4 & 30 & 3 & 3 & 31 & 18 & 16\phantom{${}^*$} & 0.29 \\
me & lh & 1,4 & 29 & 3 & 3 & 30 & 17 & 16\phantom{${}^*$} & 0.27 \\
se & ll & 2,4 & 14 & 2 & 2 & 15 & 7 & 10\phantom{${}^*$} & 0.2 \\
me & ll & 2,4 & 14 & 2 & 2 & 15 & 7 & 10\phantom{${}^*$} & 0.2 \\
se & la & 2,4 & 14 & 2 & 2 & 15 & 7 & 10\phantom{${}^*$} & 0.2 \\
me & la & 2,4 & 14 & 2 & 2 & 15 & 7 & 10\phantom{${}^*$} & 0.21 \\
se & lh & 2,4 & 29 & 3 & 3 & 30 & 17 & 16\phantom{${}^*$} & 0.29 \\
me & lh & 2,4 & 30 & 3 & 3 & 31 & 18 & 16\phantom{${}^*$} & 0.27 \\
se & ll & 3,4 & 30 & 3 & 3 & 31 & 18 & 16\phantom{${}^*$} & 0.3 \\
me & ll & 3,4 & 29 & 3 & 3 & 30 & 17 & 16\phantom{${}^*$} & 0.33 \\
se & la & 3,4 & 29 & 3 & 3 & 30 & 17 & 16\phantom{${}^*$} & 0.3 \\
me & la & 3,4 & 29 & 3 & 3 & 30 & 17 & 16\phantom{${}^*$} & 0.31 \\
se & lh & 3,4 & 30 & 3 & 3 & 31 & 18 & 16\phantom{${}^*$} & 0.33 \\
me & lh & 3,4 & 30 & 3 & 3 & 31 & 18 & 16\phantom{${}^*$} & 0.32 \\
se & ll & 1,2,3,4 & 29 & 3 & 3 & 30 & 17 & 16\phantom{${}^*$} & 0.32 \\
me & ll & 1,2,3,4 & 29 & 3 & 3 & 30 & 17 & 16\phantom{${}^*$} & 0.24 \\
se & la & 1,2,3,4 & 29 & 3 & 3 & 30 & 17 & 16\phantom{${}^*$} & 0.34 \\
me & la & 1,2,3,4 & 30 & 3 & 3 & 31 & 18 & 16\phantom{${}^*$} & 0.24 \\
se & lh & 1,2,3,4 & 30 & 3 & 3 & 31 & 18 & 16\phantom{${}^*$} & 0.28 \\
me & lh & 1,2,3,4 & 30 & 3 & 3 & 31 & 18 & 16\phantom{${}^*$} & 0.3 \\
se & ll & 5 & 125 & 6 & 6 & 126 & 53 & 32${}^*$ & 1.96 \\
me & ll & 5 & 125 & 6 & 6 & 126 & 53 & 32${}^*$ & 1.7 \\
se & la & 5 & 125 & 6 & 6 & 126 & 53 & 32${}^*$ & 1.92 \\
me & la & 5 & 125 & 6 & 6 & 126 & 53 & 32${}^*$ & 1.8 \\
se & lh & 5 & 125 & 6 & 6 & 126 & -- & --\phantom{${}^*$} & 2.42 \\
me & lh & 5 & 125 & 6 & 6 & 126 & -- & --\phantom{${}^*$} & 2.09 \\
se & ll & 1,5 & 125 & 6 & 6 & 126 & 53 & 32${}^*$ & 1.78 \\
me & ll & 1,5 & 125 & 6 & 6 & 126 & 53 & 32${}^*$ & 1.52 \\
se & la & 1,5 & 125 & 6 & 6 & 126 & 53 & 32${}^*$ & 1.5 \\
me & la & 1,5 & 125 & 6 & 6 & 126 & 53 & 32${}^*$ & 1.8 \\
se & lh & 1,5 & 125 & 6 & 6 & 126 & -- & --\phantom{${}^*$} & 2.55 \\
me & lh & 1,5 & 125 & 6 & 6 & 126 & -- & --\phantom{${}^*$} & 2.44 \\
se & ll & 2,5 & 125 & 6 & 6 & 126 & -- & --\phantom{${}^*$} & 2.27 \\
me & ll & 2,5 & 125 & 6 & 6 & 126 & -- & --\phantom{${}^*$} & 2.76 \\
se & la & 2,5 & 125 & 6 & 6 & 126 & -- & --\phantom{${}^*$} & 2.49 \\
me & la & 2,5 & 125 & 6 & 6 & 126 & -- & --\phantom{${}^*$} & 2.37 \\
se & lh & 2,5 & 125 & 6 & 6 & 126 & -- & --\phantom{${}^*$} & 2.62 \\
me & lh & 2,5 & 125 & 6 & 6 & 126 & -- & --\phantom{${}^*$} & 2.48 \\
se & ll & 3,5 & 125 & 6 & 6 & 126 & 53 & 32${}^*$ & 1.94 \\
me & ll & 3,5 & 125 & 6 & 6 & 126 & 53 & 32${}^*$ & 1.97 \\
se & la & 3,5 & 125 & 6 & 6 & 126 & 53 & 32${}^*$ & 1.82 \\
me & la & 3,5 & 125 & 6 & 6 & 126 & 53 & 32${}^*$ & 2.18 \\
se & lh & 3,5 & 125 & 6 & 6 & 126 & -- & --\phantom{${}^*$} & 2.46 \\
me & lh & 3,5 & 125 & 6 & 6 & 126 & -- & --\phantom{${}^*$} & 2.52 \\
se & ll & 4,5 & 125 & 6 & 6 & 126 & -- & --\phantom{${}^*$} & 2.87 \\
me & ll & 4,5 & 125 & 6 & 6 & 126 & -- & --\phantom{${}^*$} & 2.8 \\
se & la & 4,5 & 125 & 6 & 6 & 126 & -- & --\phantom{${}^*$} & 2.71 \\
me & la & 4,5 & 125 & 6 & 6 & 126 & -- & --\phantom{${}^*$} & 2.72 \\
se & lh & 4,5 & 125 & 6 & 6 & 126 & -- & --\phantom{${}^*$} & 2.95 \\
me & lh & 4,5 & 125 & 6 & 6 & 126 & -- & --\phantom{${}^*$} & 2.8 \\
se & ll & 1,3,5 & 125 & 6 & 6 & 126 & 53 & 32${}^*$ & 1.83 \\
me & ll & 1,3,5 & 125 & 6 & 6 & 126 & 53 & 32${}^*$ & 1.8 \\
se & la & 1,3,5 & 125 & 6 & 6 & 126 & 53 & 32${}^*$ & 1.85 \\
me & la & 1,3,5 & 125 & 6 & 6 & 126 & 53 & 32${}^*$ & 1.78 \\
se & lh & 1,3,5 & 125 & 6 & 6 & 126 & -- & --\phantom{${}^*$} & 2.55 \\
me & lh & 1,3,5 & 125 & 6 & 6 & 126 & -- & --\phantom{${}^*$} & 2.73 \\
se & ll & 1,2,3,4,5 & 125 & 6 & 6 & 126 & -- & --\phantom{${}^*$} & 3.17 \\
me & ll & 1,2,3,4,5 & 125 & 6 & 6 & 126 & -- & --\phantom{${}^*$} & 2.98 \\
se & la & 1,2,3,4,5 & 125 & 6 & 6 & 126 & -- & --\phantom{${}^*$} & 2.96 \\
me & la & 1,2,3,4,5 & 125 & 6 & 6 & 126 & -- & --\phantom{${}^*$} & 2.83 \\
se & lh & 1,2,3,4,5 & 125 & 6 & 6 & 126 & -- & --\phantom{${}^*$} & 3.03 \\
me & lh & 1,2,3,4,5 & 125 & 6 & 6 & 126 & -- & --\phantom{${}^*$} & 2.88 \\
se & ll & 6 & 145 & 6 & 6 & 146 & 116 & 36\phantom{${}^*$} & 52.3 \\
me & ll & 6 & 145 & 6 & 6 & 146 & 116 & 36\phantom{${}^*$} & 49.4 \\
se & la & 6 & 145 & 6 & 6 & 146 & 116 & 36\phantom{${}^*$} & 52.6 \\
me & la & 6 & 145 & 6 & 6 & 146 & 116 & 36\phantom{${}^*$} & 49.5 \\
se & lh & 6 & 425 & 10 & 10 & 426 & -- & --\phantom{${}^*$} & 250 \\
me & lh & 6 & 425 & 10 & 10 & 426 & -- & --\phantom{${}^*$} & 265 \\
se & ll & 1,6 & 417 & 10 & 10 & 418 & -- & --\phantom{${}^*$} & 182 \\
me & ll & 1,6 & 417 & 10 & 10 & 418 & -- & --\phantom{${}^*$} & 170 \\
se & la & 1,6 & 417 & 10 & 10 & 418 & -- & --\phantom{${}^*$} & 186 \\
me & la & 1,6 & 417 & 10 & 10 & 418 & -- & --\phantom{${}^*$} & 186 \\
se & lh & 1,6 & 425 & 10 & 10 & 426 & -- & --\phantom{${}^*$} & 265 \\
me & lh & 1,6 & 425 & 10 & 10 & 426 & -- & --\phantom{${}^*$} & 233 \\
se & ll & 2,6 & 145 & 6 & 6 & 146 & 116 & 36\phantom{${}^*$} & 55 \\
me & ll & 2,6 & 145 & 6 & 6 & 146 & 116 & 36\phantom{${}^*$} & 58 \\
se & la & 2,6 & 145 & 6 & 6 & 146 & 116 & 36\phantom{${}^*$} & 55.4 \\
me & la & 2,6 & 145 & 6 & 6 & 146 & 116 & 36\phantom{${}^*$} & 51.6 \\
se & lh & 2,6 & 425 & 10 & 10 & 426 & -- & --\phantom{${}^*$} & 264 \\
me & lh & 2,6 & 425 & 10 & 10 & 426 & -- & --\phantom{${}^*$} & 261 \\
se & ll & 3,6 & 423 & 10 & 10 & 424 & -- & --\phantom{${}^*$} & 262 \\
me & ll & 3,6 & 423 & 10 & 10 & 424 & -- & --\phantom{${}^*$} & 277 \\
se & la & 3,6 & 423 & 10 & 10 & 424 & -- & --\phantom{${}^*$} & 259 \\
me & la & 3,6 & 423 & 10 & 10 & 424 & -- & --\phantom{${}^*$} & 298 \\
se & lh & 3,6 & 425 & 10 & 10 & 426 & -- & --\phantom{${}^*$} & 276 \\
me & lh & 3,6 & 425 & 10 & 10 & 426 & -- & --\phantom{${}^*$} & 268 \\
se & ll & 4,6 & 145 & 6 & 6 & 146 & 116 & 36\phantom{${}^*$} & 63.5 \\
me & ll & 4,6 & 145 & 6 & 6 & 146 & 116 & 36\phantom{${}^*$} & 67.7 \\
se & la & 4,6 & 145 & 6 & 6 & 146 & 116 & 36\phantom{${}^*$} & 60.6 \\
me & la & 4,6 & 145 & 6 & 6 & 146 & 116 & 36\phantom{${}^*$} & 65.4 \\
se & lh & 4,6 & 425 & 10 & 10 & 426 & -- & --\phantom{${}^*$} & 308 \\
me & lh & 4,6 & 425 & 10 & 10 & 426 & -- & --\phantom{${}^*$} & 311 \\
se & ll & 5,6 & 425 & 10 & 10 & 426 & -- & --\phantom{${}^*$} & 378 \\
me & ll & 5,6 & 425 & 10 & 10 & 426 & -- & --\phantom{${}^*$} & 315 \\
se & la & 5,6 & 425 & 10 & 10 & 426 & -- & --\phantom{${}^*$} & 361 \\
me & la & 5,6 & 425 & 10 & 10 & 426 & -- & --\phantom{${}^*$} & 326 \\
se & lh & 5,6 & 425 & 10 & 10 & 426 & -- & --\phantom{${}^*$} & 344 \\
me & lh & 5,6 & 425 & 10 & 10 & 426 & -- & --\phantom{${}^*$} & 302 \\
se & ll & 2,4,6 & 145 & 6 & 6 & 146 & 116 & 36\phantom{${}^*$} & 73.9 \\
me & ll & 2,4,6 & 145 & 6 & 6 & 146 & 116 & 36\phantom{${}^*$} & 68.9 \\
se & la & 2,4,6 & 145 & 6 & 6 & 146 & 116 & 36\phantom{${}^*$} & 59.2 \\
me & la & 2,4,6 & 145 & 6 & 6 & 146 & 116 & 36\phantom{${}^*$} & 69.3 \\
se & lh & 2,4,6 & 425 & 10 & 10 & 426 & -- & --\phantom{${}^*$} & 293 \\
me & lh & 2,4,6 & 425 & 10 & 10 & 426 & -- & --\phantom{${}^*$} & 300 \\
se & ll & 1,2,3,4,5,6 & 425 & 10 & 10 & 426 & -- & --\phantom{${}^*$} & 447 \\
me & ll & 1,2,3,4,5,6 & 425 & 10 & 10 & 426 & -- & --\phantom{${}^*$} & 402 \\
se & la & 1,2,3,4,5,6 & 425 & 10 & 10 & 426 & -- & --\phantom{${}^*$} & 387 \\
me & la & 1,2,3,4,5,6 & 425 & 10 & 10 & 426 & -- & --\phantom{${}^*$} & 509 \\
se & lh & 1,2,3,4,5,6 & 425 & 10 & 10 & 426 & -- & --\phantom{${}^*$} & 397 \\
me & lh & 1,2,3,4,5,6 & 425 & 10 & 10 & 426 & -- & --\phantom{${}^*$} & 547 \\
se & ll & 7 & 1683 & 20 & 20 & 1684 & -- & --\phantom{${}^*$} & 3.22e+04 \\
me & ll & 7 & 1683 & 20 & 20 & 1684 & -- & --\phantom{${}^*$} & 2.66e+04 \\
se & la & 7 & 1683 & 20 & 20 & 1684 & -- & --\phantom{${}^*$} & 5.16e+04 \\
me & la & 7 & 1683 & 20 & 20 & 1684 & -- & --\phantom{${}^*$} & 2.02e+04 \\
se & lh & 7 & 1683 & 20 & 20 & 1684 & -- & --\phantom{${}^*$} & 3.46e+04 \\
me & lh & 7 & 1683 & 20 & 20 & 1684 & -- & --\phantom{${}^*$} & 3.06e+04 \\
se & ll & 1,7 & 1683 & 20 & 20 & 1684 & -- & --\phantom{${}^*$} & 3.11e+04 \\
me & ll & 1,7 & 1683 & 20 & 20 & 1684 & -- & --\phantom{${}^*$} & 1.65e+04 \\
se & la & 1,7 & 1683 & 20 & 20 & 1684 & -- & --\phantom{${}^*$} & 1.82e+04 \\
me & la & 1,7 & 1683 & 20 & 20 & 1684 & -- & --\phantom{${}^*$} & 5.93e+04 \\
se & lh & 1,7 & 1683 & 20 & 20 & 1684 & -- & --\phantom{${}^*$} & 5.15e+04 \\
me & lh & 1,7 & 1683 & 20 & 20 & 1684 & -- & --\phantom{${}^*$} & 4.11e+04 \\
se & ll & 2,7 & 1683 & 20 & 20 & 1684 & -- & --\phantom{${}^*$} & 7.68e+04 \\
me & ll & 2,7 & 1683 & 20 & 20 & 1684 & -- & --\phantom{${}^*$} & 4.39e+04 \\
se & la & 2,7 & 1683 & 20 & 20 & 1684 & -- & --\phantom{${}^*$} & 4.43e+04 \\
me & la & 2,7 & 1683 & 20 & 20 & 1684 & -- & --\phantom{${}^*$} & 3.22e+04 \\
se & lh & 2,7 & 1683 & 20 & 20 & 1684 & -- & --\phantom{${}^*$} & 4.17e+04 \\
me & lh & 2,7 & 1683 & 20 & 20 & 1684 & -- & --\phantom{${}^*$} & 5.04e+04 \\
se & ll & 3,7 & 1683 & 20 & 20 & 1684 & -- & --\phantom{${}^*$} & 9.72e+04 \\
me & ll & 3,7 & 1683 & 20 & 20 & 1684 & -- & --\phantom{${}^*$} & 5.89e+04 \\
se & la & 3,7 & 1683 & 20 & 20 & 1684 & -- & --\phantom{${}^*$} & 1.46e+04 \\
me & la & 3,7 & 1683 & 20 & 20 & 1684 & -- & --\phantom{${}^*$} & 1.38e+05 \\
se & lh & 3,7 & 1683 & 20 & 20 & 1684 & -- & --\phantom{${}^*$} & 3.73e+04 \\
me & lh & 3,7 & 1683 & 20 & 20 & 1684 & -- & --\phantom{${}^*$} & 3.80e+04 \\
se & ll & 4,7 & 1683 & 20 & 20 & 1684 & -- & --\phantom{${}^*$} & 3.65e+04 \\
me & ll & 4,7 & 1683 & 20 & 20 & 1684 & -- & --\phantom{${}^*$} & 3.50e+04 \\
se & la & 4,7 & 1683 & 20 & 20 & 1684 & -- & --\phantom{${}^*$} & 6.28e+04 \\
me & la & 4,7 & 1683 & 20 & 20 & 1684 & -- & --\phantom{${}^*$} & 3.75e+04 \\
se & lh & 4,7 & 1683 & 20 & 20 & 1684 & -- & --\phantom{${}^*$} & 3.69e+04 \\
me & lh & 4,7 & 1683 & 20 & 20 & 1684 & -- & --\phantom{${}^*$} & 3.61e+04 \\
se & ll & 5,7 & 1683 & 20 & 20 & 1684 & -- & --\phantom{${}^*$} & 3.48e+04 \\
me & ll & 5,7 & 1683 & 20 & 20 & 1684 & -- & --\phantom{${}^*$} & 1.09e+05 \\
se & la & 5,7 & 1683 & 20 & 20 & 1684 & -- & --\phantom{${}^*$} & 4.17e+04 \\
me & la & 5,7 & 1683 & 20 & 20 & 1684 & -- & --\phantom{${}^*$} & 4.87e+04 \\
se & lh & 5,7 & 1683 & 20 & 20 & 1684 & -- & --\phantom{${}^*$} & 4.39e+04 \\
me & lh & 5,7 & 1683 & 20 & 20 & 1684 & -- & --\phantom{${}^*$} & 4.61e+04 \\
se & ll & 6,7 & 1683 & 20 & 20 & 1684 & -- & --\phantom{${}^*$} & 4.67e+04 \\
me & ll & 6,7 & 1683 & 20 & 20 & 1684 & -- & --\phantom{${}^*$} & 4.60e+04 \\
se & la & 6,7 & 1683 & 20 & 20 & 1684 & -- & --\phantom{${}^*$} & 4.97e+04 \\
me & la & 6,7 & 1683 & 20 & 20 & 1684 & -- & --\phantom{${}^*$} & 4.76e+04 \\
se & lh & 6,7 & 1683 & 20 & 20 & 1684 & -- & --\phantom{${}^*$} & 4.10e+04 \\
me & lh & 6,7 & 1683 & 20 & 20 & 1684 & -- & --\phantom{${}^*$} & 5.31e+04 \\
se & ll & 1,3,5,7 & 1683 & 20 & 20 & 1684 & -- & --\phantom{${}^*$} & 2.63e+04 \\
me & ll & 1,3,5,7 & 1683 & 20 & 20 & 1684 & -- & --\phantom{${}^*$} & 8.89e+04 \\
se & la & 1,3,5,7 & 1683 & 20 & 20 & 1684 & -- & --\phantom{${}^*$} & 1.96e+04 \\
me & la & 1,3,5,7 & 1683 & 20 & 20 & 1684 & -- & --\phantom{${}^*$} & 2.44e+04 \\
se & lh & 1,3,5,7 & 1683 & 20 & 20 & 1684 & -- & --\phantom{${}^*$} & 3.90e+04 \\
me & lh & 1,3,5,7 & 1683 & 20 & 20 & 1684 & -- & --\phantom{${}^*$} & 3.89e+04 \\
se & ll & 1,2,3,4,5,6,7 & 1683 & 20 & 20 & 1684 & -- & --\phantom{${}^*$} & 4.67e+04 \\
me & ll & 1,2,3,4,5,6,7 & 1683 & 20 & 20 & 1684 & -- & --\phantom{${}^*$} & 4.99e+04 \\
se & la & 1,2,3,4,5,6,7 & 1683 & 20 & 20 & 1684 & -- & --\phantom{${}^*$} & 3.86e+04 \\
me & la & 1,2,3,4,5,6,7 & 1683 & 20 & 20 & 1684 & -- & --\phantom{${}^*$} & 4.95e+04 \\
se & lh & 1,2,3,4,5,6,7 & 1683 & 20 & 20 & 1684 & -- & --\phantom{${}^*$} & 5.53e+04 \\
me & lh & 1,2,3,4,5,6,7 & 1683 & 20 & 20 & 1684 & -- & --\phantom{${}^*$} & 6.12e+04

\end{longtable}

\end{document}